\documentclass[10pt]{amsart}
\usepackage{amscd}
\usepackage{amssymb}
\newtheorem{Theorem}{Theorem}[section]
\newtheorem{Lemma}[Theorem]{Lemma}
\newtheorem{Proposition}[Theorem]{Proposition}
\newtheorem{Corollary}[Theorem]{Corollary}
\newtheorem{Definition}[Theorem]{Definition}
\newtheorem{Example}[Theorem]{Example}
\newtheorem{Remark}[Theorem]{Remark}

\def\ot{\otimes}
\def\ep{\varepsilon}

\def\cok{\operatorname{cok}}

\def\coeq{\operatorname{coeq}}

\def\Cor{\operatorname{Cor}}
\def\Rad{\operatorname{Rad}}
\def\gr{\operatorname{gr}}
\def\Hom{\operatorname{Hom}}

\def\diag{\operatorname{diag}}
\def\Aut{\operatorname{Aut}}
\def\Coalg{\operatorname{Coalg}}
\def\Alg{\operatorname{Alg}}
\def\Vect{\operatorname{Vect}}

\def\res{\operatorname{res}}

\def\co-Lift{\operatorname{co-Lift}}
\def\and{\operatorname{and}}

\def\im{\operatorname{im}}

\def\res{\operatorname{res}}
\def\m{\operatorname{m}}

\def\YD{\operatorname{YD}}
\def\cal{\mathcal}
\def\incl{\operatorname{incl}}
\def\proj{\operatorname{proj}}
\def\ad{\operatorname{ad}}
\def\adj{\operatorname{adj}}
\def\coadj{\operatorname{coadj}}
\def\ad{\operatorname{ad}}
\def\coeq{\operatorname{coeq}}

\def\im{\operatorname{im}}
\def\genst#1#2{\left\langle\left. #1 \right| #2 \right\rangle}
\def\setst#1#2{\left\{\left. #1 \right| #2 \right\}}
\def\gen#1{\left\langle #1 \right\rangle}
\def\set#1{\left\{ #1 \right\}}
\def\D{\mathcal{D}}
\def\m{\mathrm{m}}
\def\id{\operatorname{id}}

\title{Pointed Hopf algebras as cocycle deformations}

\author{L. Grunenfelder and M. Mastnak}

\address{Department of Mathematics, The University of British Columbia,
Vancouver, BC V6T 1Z2, Canada}

\email{luzius@math.ubc.ca, luzius@mathstat.dal.ca}

\address{Department of Mathematics and C.S., Saint Mary's University, Halifax,
NS B3H 3C3, Canada}

\email{mmastnak@cs.smu.ca}

\thanks{Research supported in part by NSERC}

\begin{document}

\tolerance=1000

\begin{abstract} We show that all finite dimensional pointed Hopf
algebras with the same diagram in the classification scheme of
Andruskiewitsch and Schneider are cocycle deformations of each
other. This is done by giving first a suitable characterization of
such Hopf algebras, which allows for the application of results by
Masuoka about Morita-Takeuchi equivalence and by Schauenburg about
Hopf Galois extensions. We also outline a method to describe the
deforming cocycles involved using the exponential map and its
q-analogue.
\end{abstract}

\maketitle

\setcounter{section}{-1}

\section{Introduction}\label{s0}

Finite dimensional pointed Hopf algebras over an algebraically
closed field of characteristic zero, particularly when the group
of points is abelian, have been studied quite extensively with
various methods in \cite{AS, BDG, Gr1, Mu}. The most far reaching
results as yet in this area have been obtained in \cite{AS}, where
a large class of such Hopf algebras are classified. In the present
paper we show, among other things, that all Hopf algebras in
this class can be obtained by cocycle deformations from Radford
biproducts of the form $B(V)\# kG$, where $B(V)$ is the Nichols
algebra of the Yetter-Drinfeld $kG$-module $V$.

After Kaplansky's tenth conjecture concerning the finiteness of
the set of isomorphism classes of Hopf algebras of a given finite
dimension had been refuted, a weakened version of that conjecture
has been proposed \cite{Ma, Di}, namely that there are only
finitely many quasi-isomorphism classes of Hopf algebras of a
given (odd) dimension. Two Hopf algebras are said to be
quasi-isomorphic if they have equivalent comodule categories. In
even dimensions the weakened conjecture has been disproved in
\cite{EG, Gra}, where it is shown that there are infinitely many
isomorphism classes of Hopf algebras of dimension 32. Our results
confirm the conjecture for a large class of pointed Hopf algebras
of odd dimension. If $H$ is a Hopf algebra with coradical a Hopf
subalgebra then the graded coalgebra $\gr^cH$ associated with the
coradical filtration is a graded Hopf algebra and its elements of
positive degree form the radical. If the radical of $H$ is a Hopf
ideal then the graded algebra associated with the radical
filtration is a graded Hopf algebra with $\Cor (\gr^r H)\cong
H/\Rad H$. In either case we have $\gr H\cong R\# H_0$, where
$H_0$ is the degree zero part and $R$ is the braided Hopf algebra
of coinvariants or invariants, respectively.

The Nichols algebra $B(V)$ of a crossed $kG$-module $V$ is a
connected graded braided Hopf algebra. The Radford biproduct
$H(V)=B(V)\# kG$ is an ordinary graded Hopf algebra with coradical
$kG$ and the elements of positive degree form a Hopf ideal (the
graded radical). A lifting of $H(V)$ is a pointed Hopf algebra $H$
for which $\gr^cH\cong H(V)$. Such liftings are obtained by
deforming the multiplication of $H(V)$. The lifting problem for
$V$ asks for the classification of all liftings of $H(V)$. This
problem, together with the characterization of $B(V)$ and $H(V)$,
have been solved by Andruskiewitsch and Schneider in \cite{AS} for
a large class of crossed $kG$-modules of finite Cartan type. It
allows them to classify all finite dimensional pointed Hopf
algebras $A$ for which the order of the abelian group of points
has no prime factors $< 11$. In this paper we find a description
of these lifted Hopf algebras, which is suitable for the
application of a result of Masuoka about Morita-Takeuchi
equivalence \cite{Ma} and of Schauenburg about Hopf Galois
extensions \cite{Sch}, to prove that all liftings of a given
$H(V)$ in this class are cocycle deformations of each other. As a
result we see here that in the class of finite dimensional pointed
Hopf algebras classified by Andruskiewitsch an Schneider \cite{AS}
all Hopf algebras $H$ with isomorphic associated graded Hopf
algebra $\gr^cH$ are monoidally Morita-Takeuchi equivalent, and
therefore cocycle deformations of each other. For some special
cases such results have been obtained in \cite{Ma, Di, BDR}.

In Section \ref{s1} we give a short review of braided spaces,
braided Hopf algebras, Nichols algebras and bosonization (among other things, we prove
Theorem \ref{quantum-symmetrizer} in which we give a new way of looking at the quantum symmetrizer).
Most of this is done in preparation for a useful characterization of the liftings for a
large class of crossed modules over finite abelian groups in
Section \ref{s2}. With this characterization it is then possible in Section \ref{s3}
to prove that liftings are Morita-Takeuchi equivalent (Theorem \ref{main4}) by using Masuoka's pushout construction, and that they
are cocycle deformations of each other (Corollary \ref{main5}) by a result of
Schauenburg.  Cocycle deformations as well as their relation to Hochschild
cohomology are discussed in Section \ref{s4}. See \cite{GM1} and \cite{Gr2} for a
more detailed account of the cohomological aspects. In particular, a complete description is given of all relevant cocycles for the quantum linear case, where exponential map and its $q$-analogue play a prominent role.  In Subsection \ref{general} we outline the main ideas of
\cite{GM2}, where we, among other things, use the exponential map, its $q$-analogue, and the theory of Singer extensions, to explicitly describe the
deforming cocycles for all Anruskiewitsch-Schneider Hopf algebras.  Explicit examples of deforming cocycles are presented in Section \ref{s5}.  We recommend to the reader to start by looking at Example \ref{example1} where the theory is illustrated in great detail on the smallest nontrivial example.

\subsection*{Acknowledgements} We thank A. Masuoka for correcting an error after reading an early
version \cite{GM1} of a part of this work.  The second author is
very grateful for useful conversations with M.~Beattie.

\section{Braided Hopf algebras and the (bi-)crossproduct}\label{s1}

Our main interest in this paper are pointed Hopf
algebras. We also briefly look at the dual notion, that is copointed
Hopf algebras.  Here we review some relevant prerequisites and facts.

\subsection{Pointed and copointed Hopf algebras} A Hopf
algebra $H$ is pointed if its coradical $\Cor H$ is equal to the
group algebra of the group of points $G(H)$. In this case the
coradical filtration is an ascending Hopf algebra filtration and
the associated graded Hopf algebra $\gr^cH$ has the obvious
injection $\kappa^c\colon kG\to\gr^cH$ and projection $\pi^c\colon
\gr^cH\to kG$ such that $\pi^c\kappa^c =1$.

We say that $H$ is copointed if its radical $\Rad H$ is a Hopf
ideal and $H/{\Rad H}$ is a group algebra $kG$. Here the radical
filtration is an descending Hopf algebra filtration and again the
associated graded Hopf algebra $\gr^rH$ has the obvious projection
$\pi^r\colon \gr^rH\to kG$ and an injection $\kappa^r\colon
kG\to\gr^rH$ such that $\pi^r\kappa^r =1$.

In both cases above $\gr H$ is graded, pointed and copointed, and
by \cite{Ra1} $\gr H\cong A\# kG$, where $A=\setst{ x\in
H}{(\pi\ot 1)\Delta (x) =1\ot x}$ is the graded connected braided
Hopf algebra of coinvariants.

\begin{Lemma}\label{l17} If the Hopf algebra $H$ is pointed
and copointed then $\Cor H\cong H/{\Rad H}$ and $H$ is a Hopf
algebra with a projection. Moreover, $R\# kG\cong H$, where
$R=\setst{x\in \gr H}{(p\ot 1)\Delta (x)=1\ot x}$ is the connected
braided Hopf algebra of coinvariants of $H$. This is the case in
particular for $\gr^cH$ and for $\gr^rH$ when $H$ is pointed or
copointed, respectively.
\end{Lemma}

\begin{proof} A surjective coalgebra map $\eta\colon C\to D$, where
$D=\Cor (D)$, maps $\Cor C$ onto $\Cor D$ \cite{Mo}. Thus, the
composite $\Cor H\to H\to H/{\Rad H}$ is a bijection. The
isomorphism is that of \cite{Ra1}.
\end{proof}

\subsection{Braidings}
A braided monoidal category $\mathcal V$ is a monoidal category
together with a natural morphism $c\colon V\ot W\to W\ot V$ such
that
\begin{enumerate}
\item $c_{k,V}=\tau =c_{V,k}$,
\item $c_{U\ot V,W}=(c_{U,W}\ot 1)(1\ot c_{V,W})$,
\item $c_{U,V\ot W}=(1\ot c_{U,W})(c_{U,V}\ot 1)$,
\item $c(f\ot g)=(g\ot f)c$.
\end{enumerate}
Braided algebras, braided coalgebras and braided Hopf algebras are
now defined with this tensor product and braiding in mind. The
compatibility condition $\Delta\m =(m\ot m)(1\ot c\ot
1)(\Delta\ot\Delta )$ between multiplication and comultiplication
in a braided Hopf algebra $A$ involves the braiding $c\colon A\ot
A\to A\ot A$, so that the diagram
$$\begin{CD}
A\ot A @> (1\ot c\ot 1)(\Delta\ot\Delta )>> A\ot A\ot A\ot A \\
@V m VV  @V m\ot m VV \\
A @> \Delta >>  A\ot A
\end{CD}$$
commutes, i.e. multiplication and unit are morphisms of braided
coalgebras or, equivalently, comultiplication and counit are maps
of braided algebras.

In the present paper $\mathcal V$ is a category of braided vector
spaces, where the braidings $c\colon V\ot V\to V\ot V$ are of
finite abelian group type, so that $c\colon V\ot V\to V\ot V$ is
of diagonal type.

\subsection{Primitives and indecomposables}

The vector space of primitives
$$P(A) =\setst{ y\in A}{\Delta (y)=y\ot 1+1\ot y}\cong\ker
(\tilde{\Delta})\colon A/k\to A/k\ot A/k$$ of a braided Hopf
algebra $A$ is a braided vector space, since $\Delta$ is  a map in
$\mathcal V$. The c-bracket map $[-,-]_c=m(1\ot 1-c)\colon A\ot
A\to A$ restricted to $P(A)$ satisfies $\Delta [x,y]_c=[x,y]_c\ot
1+(1-c^2)x\ot y +1\ot [x,y]_c$; in particular $[x,y]_c\in P(A)$ if
and only $c^2(x\ot y)=x\ot y$. Moreover, if $x\in P(A)$ and
$c(x\ot x)=qx\ot x$ then $\Delta x^n=\sum_{i+j=n}{n\choose
i}_qx^i\ot x^j$, where ${n\choose
i}_q={\frac{n_q!}{i_q!(n-i)_q!}}$ are the $q$-binomial
coefficients (the Gauss polynomials) for $q$, $m_q!=1_q2_q...m_q$
with $j_q=1+q+...+q^{j-1}$ if $j>0$ and $0_q!=1$. If $q=1$ then
$j_q=j$ and we have the ordinary binomial coefficients, otherwise
$j_q={\frac{1-q^j}{1-q}}$. In particular, if $q$ has order $n$
then ${n\choose i}_q=0$ for $0<i<n$, and hence $x^n\in P(A)$.

\begin{Lemma} Let $\set{x_i}$ be a basis of $P(A)$ such that
$c(x_i\ot x_j)=q_{ji}x_j\ot x_i$. If $q_{ji}q_{ij}q_{ii}^{r-1}=1$ then
$\ad x_i^r(x_j)$ is primitive.
\end{Lemma}

\begin{proof} See for example \cite[Appendix 1]{AS1}.
\end{proof}

The vector space of indecomposables
$$Q(A)=JA/JA^2=\cok (\tilde m\colon JA\ot JA\to JA),$$
where $JA=\ker (\epsilon )$, is a braided vector space as well.
The c-cobracket map $\delta_c=(1\ot 1-c)\Delta\colon A\to A\ot A$
restricts to $JA$, since
$$\Delta (\bar a)= \bar a\ot 1+1\ot\bar a
+\sum \bar a_i\ot\bar b_i$$
and hence
$$\delta (\bar a)=\sum (\bar a_i\ot\bar b_i -c(\bar a_i\ot\bar b_i))$$
is in $JA\ot JA$ for every $\bar a=a-\epsilon (a)\in JA$. Moreover,
$$\delta (\bar a\bar b)-(\bar a\ot\bar b -c^2(\bar a\ot\bar b))$$ is
in $JA^2\ot JA + JA\ot JA^2$. In particular, if $c^2(\bar a\ot\bar
b)=\bar a\ot\bar b$, then $\delta (\bar a\ot\bar b)\in JA^2\ot
JA+JA\ot JA^2$.

\subsection{The free and the cofree graded braided Hopf algebras}

The forgetful functor $U\colon \Alg_c\to \mathcal V_c$ has a
left-adjoint $\cal A\colon \mathcal V_c\to \Alg_c$ and the
forgetful functor $U\colon \Coalg_c\to\mathcal V_c $ has a
right-adjoint $\cal C\colon \mathcal V_c\to \Coalg_c$, the free
braided graded algebra functor and the cofree graded braided
coalgebra functor, respectively. Moreover, there is a natural
transformation $\cal S\colon \cal A\to\cal C$, the shuffle map or
quantum symmetrizer. They can be described as follows.

If $(V, \mu , \delta )$ is a braided vector space then the tensor
powers $T_0(V)=k$, $T_{n+1}(V)=V\ot T_n(V)$ are braided vector
spaces as well and so is $T(V)=\oplus_nT_n(V)$. The ordinary
tensor algebra structure makes $T(V)$ the free connected graded
braided algebra, and the ordinary tensor coalgebra structure makes
it the cofree  connected graded braided coalgebra.

By the universal property of the graded braided tensor algebra
$T(V)$ the linear map $\Delta_1=\incl\diag\colon V\to T(V)\ot
T(V)$, $\Delta_1 (v)=v\ot 1+1\ot v$, induces the c-shuffle
comultiplication $\Delta_{\cal A}\colon T(V)\to T(V)\ot T(V)$,
which is a homomorphism of braided algebras, so that $\Delta_{\cal
A}m=(m\ot m)(1\ot c\ot 1)(\Delta_{\cal A}\ot\Delta_{\cal A})$.
Moreover, the linear map $s_1\colon V\to T(V)$, $s_1(v)=-v$,
extends uniquely to a c-antipode $s_{\cal A}\colon T(V)\to T(V)$,
such that $s_{\cal A}m=m(s_{\cal A}\ot s_{\cal A})c$,
$\Delta_{\cal A}s_{\cal A}=c(s_{\cal A}\ot s_{\cal A})\Delta_{\cal
A}$ and $m(1\ot s_{\cal A})\Delta_{\cal A}=\iota\epsilon
=m(s_{\cal A}\ot 1)\Delta_{\cal A}$, thus making $\cal A(V)=(T(V),
m,\Delta_{\cal A},s_{\cal A})$ the free connected graded braided
Hopf algebra. This defines a functor $\cal A\colon \mathcal V_c\to
Hopf_c$, left-adjoint to the space of primitives functor $P\colon
Hopf_c\to\mathcal V_c$.

On the other hand, by the universal property of the cofree
connected graded braided coalgebra $T(V)$ there is a unique
c-shuffle multiplication $m_{\cal C}\colon T(V)\ot T(V)\to T(V)$,
which is the homomorphism of braided coalgebras induced by the
linear map $m_1=+\proj\colon T(V)\ot T(V)\to V$, so that $\Delta
m_{\cal C}=(m_{\cal C}\ot m_{\cal C})(1\ot c\ot 1)(\Delta\ot\Delta
)$. The linear map $s_1=-\proj\colon  T(V)\to V$ induces uniquely
a c-antipode $s_{\cal C}\colon T(V)\to T(V)$, such that $s_{\cal
C}m_{\cal C}=m_{\cal C}(s_{\cal C}\ot s_{\cal C})c$, $\Delta
s_{\cal C}=c(s_{\cal C}\ot s_{\cal C})\Delta$ and $m_{\cal C}(1\ot
s_{\cal C})\Delta =\iota\epsilon =m_{\cal C}(s_{\cal C}\ot
1)\Delta$, making $\cal C(V)=(T(V), \Delta , m_{\cal C}, s_{\cal
C})$ the cofree connected graded braided Hopf algebra. The functor
$\cal C \colon \mathcal V_c\to Hopf_c$ is right-adjoint to the
space of indecomposables functor $Q\colon Hopf_c\to\mathcal V_c$.

\begin{Theorem}\label{quantum-symmetrizer} There is a natural transformation
$$\mathcal S\colon \mathcal A\to \mathcal C,$$
the quantum symmetrizer, such that
$$\mathcal B(V)\cong\mathcal A(V)/\ker (\mathcal S)
\cong \im\mathcal S \subset\mathcal C(V)$$
is the Nichols algebra of $V$ and $Q\mathcal B(V)\cong V\cong
P\mathcal B(V)$. In particular, the Hopf ideal of $\mathcal A(V)$
generated by the primitives of degree $\ge 2$ is contained in
$\ker\mathcal S$.
\end{Theorem}

\begin{proof} The adjunctions just described provide natural isomorphisms
$$\mathcal{V}_c(Q\cal A(V),W)\cong Hopf_c(\cal A(V),
\cal C(W))\cong \cal V_c(V,P\cal C(W)).$$
By construction we also have
$$Q\cal A (V)\cong V \quad , \quad W\cong P\cal C (W),$$
The resulting natural isomorphism
$$\theta_{V,W}\colon \mathcal V_c(V,W)\to Hopf_c(\cal A(V),\cal C(W))$$
sends the identity morphism of $V$ to the quantum symmetrizer
$$\cal S=\theta_{V,V}(1_V)\colon \cal A(V)\to\cal C(V).$$
The image of $\cal S$ is the
Nichols algebra
$$\cal B(V)\cong\cal A(V)/\ker (\cal S)\cong\im\cal S\subset\cal C(V)$$
and $Q\cal B(V)\cong V\cong P\cal B(V)$. Moreover, since $\mathcal
S$ is graded, it follows that $\mathcal S(y)=0$ for every
primitive $y\in\mathcal A(V)$ of degree $\ge 2$.
\end{proof}

An explicit description of the quantum symmetrizer can be obtained
directly in term of the action of the braid groups $B_n$ on the tensor
powers $V^{\ot n}$.
\subsection{Crossed modules} A prime example of a
braided monoidal category is the category of crossed
$H$-modules $\YD^H_H$ for a Hopf algebra $H$. A crossed $H$-module
or a Yetter-Drinfeld $H$-module, $(V,\mu ,\delta )$ is a vector
space $V$ with a $H$-module structure $\mu \colon H\ot V\to V$,
$\mu (h\ot v)=hv$, and a $H$-comodule structure $\delta\colon V\to
H\ot V$, $\delta (v)=v_{-1}\ot v_0$, such that $h\delta
(v)=h_1v_{-1}\ot h_2v_0=(h_1v)_{-1}h_2\ot (h_1v)_0$, or
$$(m\ot\mu )(1\ot\tau\ot 1)(\Delta\ot\delta )=
(m\ot 1)(1\ot\tau )(\delta\mu\ot 1)(1\ot\tau )(\Delta\ot 1),$$
i.e such that the diagram

$$\begin{CD}
H\ot V @>\Delta\ot\delta >> H\ot H\ot H\ot V @>1\ot\tau\ot 1 >> H\ot H\ot
H\ot V \\
@V(1\ot\tau )(\Delta\ot 1)VV  @.  @Vm\ot\mu VV \\
H\ot V\ot H @>\delta\mu\ot 1 >> H\ot V\ot H @>(m\ot 1)(1\ot\tau )>> H\ot V
\end{CD}$$
commutes. This is the case in particular when
$\delta (hv)=h_1v_{-1}s(h_3)\ot h_2v_0$, i.e: when the diagram

$$\begin{CD}
H\ot V @> (1\ot\phi\ot 1)(\Delta\ot\delta ) >> H\ot H\ot H\ot V\\
@V\mu VV  @V m\ot\mu VV \\
V @> \delta >> H\ot V
\end{CD}$$
commutes, where $\phi = (m\ot 1)(1\ot s\ot 1)(1\ot\tau\Delta
)\tau\colon H\ot H\to H\ot H$, $\phi (g\ot h)=hs(g_2)\ot g_1$. These
braided $H$-modules with the obvious homomorphisms form a braided
monoidal category, with the ordinary tensor product of vector
spaces together with diagonal action and diagonal coaction. The
braiding, given by $c(v\ot w)=v_{-1}(w)\ot v_0$,
$$c=(\mu\ot 1)(1\ot\tau )(\delta\ot 1)\colon V\ot W\to W\ot V,$$
clearly satisfies the braiding conditions. The crossed $H$-module
$(k, \mu =\ep\ot 1, \delta =\iota\ot 1)$ acts as a unit for the
tensor. Moreover, $(H,\adj ,\Delta )$ and $(H, m, \coadj )$ are
crossed $H$-modules, where $\adj (h\ot h')=h_1h'S(h_2)$ and
$\coadj (h)=h_1S(h_3)\ot h_2$.

If $H=kG$ for some finite abelian group and $V$ is finite
dimensional, then the action of $G$ is diagonalizable, so that
$$V\cong\oplus_{g\in G, \chi\in\hat G}V_g^{\chi},$$
where $V_g^{\chi}=V_g\cap V^{\chi}$ with $V_g=\{ v\in V|\delta
(v)=g\ot v\}$ and $V^{\chi}=\{ v\in V| gv=\chi (g)v, \forall g\in
G\}$.

\subsection{The pushout construction for bi-cross products}
Recall Masuoka's pushout construction for Hopf algebras \cite{Ma},
\cite{Gr1}. If $A$ is a Hopf algebra then $\Alg (A,k)$  is a group
under convolution which acts on $A$ by conjugation as Hopf algebra
automorphisms.

\begin{Lemma} For every Hopf algebra $A$ the group $\Alg (A,k)$ acts on
$A$
by `conjugation' as Hopf algebra automorphisms
$$\rho\colon \Alg(A,k)\to\Aut_{Hopf}(A)^{op},$$
where $\rho_f=f*1*fs$, i.e: $\rho_f(x)=f(x_1)x_2f(sx_3)$. The
image of $\rho$ is a normal subgroup of $\Aut_{Hopf}(A)$.
\end{Lemma}

\begin{proof} It is easy to verify that $\rho_f$ is an Hopf algebra
map. The definition of $\rho$ shows
that $\rho_{f_1*f_2}=(f_1*f_2)*1*(f_2s*f_1s)=\rho_{f_2}\rho_{f_1}$
and, since $f*fs=\ep =fs*f$, it follows that
$\rho_f\rho_{fs}=1=\rho_{fs}\rho_f$, so that $\rho_f$ is a Hopf
algebra automorphism. If $\phi\in\Aut_{Hopf}(A)$ and
$f\in\Alg(A,k)$ then $\phi^{-1}\rho_f\phi
=f\phi*1*fs\phi=\rho_{f\phi}$, hence the image of $\rho$ is a
normal subgroup.
\end{proof}

Two Hopf ideals $I$ and $J$ of $A$ are said to be conjugate if
$J=\rho_f(I)=f*I*fs$ for some $f\in\Alg (A,k)$. If $x\in P_{1,g}$
is a $(1,g)$-primitive then
$$\rho_f(x)=f(x)+f(g)x+f(g)gfs(x)=f(g)x+f(x)(g-1).$$

\begin{Theorem}\cite[Theorem 2]{Ma}\cite[Theorem 3.4]{BDR}\label{po}
Let $A'$ be a Hopf subalgebra of $A$. If the Hopf ideals $I$ and
$J$ of $A'$ are conjugate and $A/(f*I)\ne 0$ then the quotient
Hopf algebras $A/(I)$ and $A/(J)$ by the Hopf ideals in $A$
generated by $I$ and $J$ are monoidally Morita-Takeuchi
equivalent, i.e: there exists a $k$-linear monoidal equivalence
between their (left) comodule categories.
\end{Theorem}

\begin{proof} Masuoka's result \cite[Theorem 2]{Ma}, that there is a
$(A/(I),A/(J))$-biGalois object, namely $A/(f*I)$, holds , provided that
$A/(f*I)\ne 0$ (\cite{BDR}, Theorem 3.4),
and we
can invoke \cite[Corollary 5.7]{Sch}, to see that $A/(I)$ and
$A/(J)$ are Morita-Takeuchi equivalent.
\end{proof}

Observe, as Masuoka did \cite{Ma}, that the commutative square
$$\begin{CD}
A' @>>> A \\
@VVV @VVV \\
B/I @>>> A/(I)
\end{CD}$$
is a pushout of Hopf algebras.

If $R$ is a braided Hopf algebra in the braided category of
crossed $H$-modules then the bi-cross product $R\# H$ is an
ordinary Hopf algebra with multiplication
$$(x\# h)(x'\# h') = xh_1(x')\# h_2h'$$
and comultiplication
$$\Delta (x\# h) = x_1\# (x_2)_{-1}h_1\ot (x_2)_0\# h_2.$$
The (left) action of $H$ on $R$ induces a (right) action on
$\Alg(R,k)$ by $fh(x)=f(hx)$. An algebra map $f\colon R\to k$ is
$H$-invariant if $fh=\ep (h)f$ for all $h\in H$.

\begin{Proposition}\label{invariant} Let $K$ be a
Hopf algebra in the braided category of crossed $H$-modules and
let $\Alg_H(K,k)$ be the set of $H$-invariant algebra maps. Then:
\begin{enumerate}
\item $\Alg_H(K,k)$ is a group under convolution.
\item The restriction map
$\res\colon  _H\Alg_H(K\# H,k)\to\Alg_H(K,k)$, $\res (F)=F\ot\iota$,
is an isomorphism of groups with inverse given by $\res^{-1}(f)=f\ot\ep$.
\item The image of the conjugation homomorphism
$$\Theta =\rho\res^{-1}\colon \Alg_H (K,k)\to \Aut_{Hopf}(K\# H)^{op}$$
is contained in $\widetilde{\Aut}_{Hopf}(K\# H)= \setst{
\phi\in\Aut_{Hopf}(K\# H)}{\phi_{|H}=id}$.
\end{enumerate}
\end{Proposition}

\begin{proof} The set of algebra maps $\Alg (K,k)$ may
not be a group, but since the coequalizer
$K^H=\coeq (\mu ,\ep\ot 1\colon H\ot K\to K)$ is an ordinary Hopf
algebra, $\Alg_H(K,k)\cong\Alg(K^H,k)$ is a group under
convolution. More directly, if $f, f'\in \Alg_H(K,k)$ then
\begin{eqnarray*}
f*f'(x y) & = & f\ot f'(x_1(x_2)_{-1}y_1\ot (x_2)_0y_2) =
f(x_1)f'(y_1)f(x_2)f'(y_2) \\
         & = & (f*f')(x)(f*f')(y)  \\
f*f'(h x) & = & f\ot f'(h_1x_1\ot h_2x_2)=f(h_1x_1)f'(h_2x_2)=\ep
(h)(f*f')(x),
\end{eqnarray*}
and $f*fs  = \ep =fs*f,$ so that $\Alg_H(K,k)$ is closed under
convolution multiplication and inversion.

For $F\in {_H{\Alg_H(K\# H,k)}}$ the map $\res(F)\colon K\to k$ is
in fact a $H$-invariant algebra map, since
\begin{eqnarray*}
\res(F)(h x)&=& F(hx\ot 1) = F((1\ot h_1)(x\ot 1)(1\ot s(h_3)) \\
&=& \ep (h)F(x\ot 1)=\ep (h)\res (F)(x)
\end{eqnarray*}
and
\begin{eqnarray*}
\res (F)(x y)=F(xy\ot 1)=F(x\ot 1)F(y\ot 1)=\res(F)(x)\res (F)(y).
\end{eqnarray*}
If $F'\in _H{\Alg_H(K\# H,k)}$ as well, then
\begin{eqnarray*}
\res(F*F')(x) &=& F\ot F'(x_1\ot (x_2)_{-1}\ot (x_2)_0\ot 1) \\
&=& F(x_1\ot 1)F'(x_2\ot 1)=\res (F)*\res (F')(x),
\end{eqnarray*}
showing that $\res$ is a group homomorphism. It is now easy to see
that $\res$ is invertible and that the inverse is as stated.

As a composite of two group homomorphisms $\Theta$ is obviously a
group homomorphism. Moreover,
$$\Theta (f)(1\ot h)=\res^{-1}(f)*1*\res^{-1}(f)s(1\ot h)=
\ep (h_1)(1\ot h_2)\ep (h_3)=1\ot h$$
for $f\in\Alg_H(K,k)$, showing that $\Theta (f)_{|H}=id$.
\end{proof}

\begin{Corollary} Let $R$ be a braided Hopf algebra in the braided
category of crossed $H$-modules and let $K$ be a braided Hopf
subalgebra. If $I$ is a Hopf ideal in $K$ and $f\in\Alg_H(K,k)$
then,
\begin{itemize}
\item $J=I\# H$ and $J_f=\Theta (f)(J)$ are Hopf ideals in $K\# H$,
\item $R\# H/(J)=R/(I)\# H$ and $R\# H/(J_f)$ are monoidally
Morita-Takeuchi equivalent, if $(R\# H)/(\res^{-1}(f)*J)\ne 0$.
\end{itemize}
\end{Corollary}

\section{Liftings over finite abelian groups}\label{s2}

In this section we give a somewhat different characterization of
the class of finite dimensional pointed Hopf algebras classified
in \cite{AS}, and show that any two such Hopf algebras with
isomorphic associated graded Hopf algebras are monoidally
Morita-Takeuchi equivalent, and therefore cocycle deformations of
each other, as we will point out in the next section.

A datum of finite Cartan type
$$\cal D =\cal D \left(G, (g_i)_{1\le i\le \theta},
(a_{ij})_{1\le i,j\le\theta}\right)$$ for a (finite) abelian group
$G$ consists of elements $g_i\in G$, $\chi_j\in \widehat{G}$ and a
Cartan matrix $(a_{ij})$ of finite type satisfying the Cartan
condition
$$q_{ij}q_{ji}=q_{ii}^{a_{ij}}$$
with $q_{ii}\ne 1$, where $q_{ij}=\chi_j(g_i)$, in particular
$q_{ii}^{a_{ij}}=q_{jj}^{a_{ji}}$ for all $1\le i,j\le \theta$. In
general, the matrix $(q_{ij})$ of a diagram of Cartan type is not
symmetric, but by \cite[Lemma 1.2]{AS} it can be reduced to the
symmetric case by twisting.

Let $\mathbf Z [I]$ be the free abelian group of rank $\theta$
with basis $\set{\alpha_1, \alpha_2,\ldots ,\alpha_{\theta}}$. The
Weyl group $W\subset \Aut (\mathbf Z[I])$ of $(a_{ij})$ is
generated by the reflections $s_i\colon \mathbf Z[I]\to \mathbf
Z[I]$, where $s_i(\alpha_j)=\alpha_j -a_{ij}\alpha_i$ for all
$i,j$. The root system of the Cartan matrix $(a_{ij})$ is $\Phi
=\cup_{i=1}^{\theta}W(\alpha_i)$ and $\Phi^+ =\Phi\cap\mathbf
Z[I]=\setst{\alpha\in\Phi}{\alpha =\sum _{i=1}^{\theta}n_i\alpha_i
, n_i\ge 0}$ is the set of positive roots relative to the basis of
simple roots $\set{\alpha_1 ,\alpha_2 ,\ldots ,\alpha_{\theta}}$.
Obviously, the number of positive roots $p=|\Phi^+|$ is at least
$\theta$. The maps $g\colon \mathbf Z[I]\to G$ and $\chi\colon
\mathbf Z[I]\to\tilde G$ given by
$g_{\alpha}=g_1^{n_1}g_2^{n_2}\ldots g_{\theta}^{n_{\theta}}$ and
$\chi_{\alpha}=\chi_1^{n_1}\chi_2^{n_2}\ldots\chi_{\theta}^{n_{\theta}}$
for $\alpha =\sum_{i=1}^{\theta}n_i\alpha_i$, respectively, are
group homomorphisms. The bilinear map $q\colon \mathbf Z
[I]\times\mathbf Z [I]\to k^{\times}$ defined by
$q_{\alpha_i\alpha_j}=q_{ij}$ can be expressed as
$q_{\alpha\beta}=\chi_{\beta}(g_{\alpha})$.

If $\cal X$ the set of connected components of the Dynkin diagram
of $\Phi$ let $\Phi_J$ be the root system of the component
$J\in\cal X$. The partition of the Dynkin diagram into connected
components corresponds to an equivalence relation on $I=\set{
1,2,\ldots ,\theta}$, where $i\sim j$ if $\alpha_i$ and $\alpha_j$
are in the same connected component.

\begin{Lemma} \cite[Lemma 2.3]{AS} Suppose that $\mathcal D$ is a
connected datum of finite Cartan type, i.e: the
Dynkin diagram of the Cartan matrix $(a_{ij})$ is connected, and such that
\begin{enumerate}
\item $q_{ii}$ has odd order, and
\item the order of $q_{ii}$ is prime to 3, if $(a_{ij})$ is of type $G_2$.
\end{enumerate}
Then there are integers $d_i\in\set{1,2,3}$ for $1\le i\le\theta$
and a $q\in k^{\times}$ of odd order $N$ such that
$$q_{ii}=q^{2d_i} \quad . \quad d_ia_{ij}=d_ja_{ji}$$
for $1\le i,j\le\theta$. If the Cartan matrix $(a_{ij})$ of
$\mathcal D$ is of type $G_2$ then the order of $q$ is prime to 3.
In particular, the $q_{ii}$ all have the same order in
$k^{\times}$, namely $N$.
\end{Lemma}

More generally, let $\mathcal D$ be a datum of finite Cartan type
in which the order $N_i$ of $q_{ii}$ is odd for all $i$, and the
order of $q_{ii}$ is prime to 3 for all $i$ in a connected
component of type $G_2$. It then follows that the order function
$N_i$ is constant, say equal to $N_J$, on each connected component
$J$. A datum satisfying these conditions will be called special
datum of finite Cartan type.

Fix a reduced decomposition of the longest element
$$w_0=s_{i_1}s_{i_2}\ldots s_{i_p}$$
of the Weyl group $W$ in terms of the simple reflections. Then
$$\set{s_{i_1}s_{i_2}\ldots s_{i_{l-1}}(\alpha_{i_l} )}_{i=1}^p$$
is a convex ordering of the positive roots.

Let $V=V(\mathcal D)$ be the crossed $kG$-module with basis $\set{
x_1,x_2,\ldots ,x_{\theta}}$, where $x_i\in V_{g_i}^{\chi_i}$ for
$1\le i\le\theta$. Then for all $1\le i\ne j\le\theta$ the elements
$ad^{1-a_{ij}}x_i(x_j)$ are primitive in the free braided
Hopf algebra $\mathcal A(V)$ (see Lemma \ref{l17} or
\cite[Appendix 1]{AS1}). If $\mathcal D$ is as in the previous
Lemma then $\chi_i^{1-a_{ij}}\chi_j\ne\ep$. This implies that
$f(u_{ij})=0$ for any braided (Hopf) subalgebra $A$ of $\mathcal
A(V)$ containing $u_{ij}=ad^{1-a_{ij}}x_i(x_j)$ and any
$G$-invariant algebra map $f\colon A\to k$. Define root vectors in
$\mathcal A(V)$ as follows by iterated braided commutators of the
elements $x_1,x_2,\ldots ,x_{\theta}$, as in Lusztig's case but
with the general braiding:
$$x_{\beta_l}=T_{i_1}T_{i_2}\ldots T_{i_{l-1}}(x_{i_l}),$$
where $T_i(x_j)=\ad_{x_i}^{-a_{ij}}(x_j)$

In the quotient Hopf algebra $R(\mathcal D)=\mathcal
A(V)/(ad^{1-a_{ij}}x_i(x_j)|1\le i\ne j\le\theta )$ define root
vectors $x_{\alpha}\in\mathcal A(V)$ for $\alpha\in\Phi^+$ by the
same iterated braided commutators of the elements $x_1, x_2,\ldots
,x_{\theta}$ as in Lusztig's case but with respect to the general
braiding. (See \cite{AS2}, and the inductive definition of root
vectors in \cite{Ri} or also \cite[Section 8.1 and Appendix]{CP}.)
Let $K(\mathcal D)$ be the subalgebra of $R(\mathcal D)$ generated
by $\setst{ x_{\alpha}^N}{\alpha\in\Phi^+}$.

\begin{Theorem}\cite[Theorem 2.6]{AS}\label{ASth2.6} Let
$\mathcal D$ be a connected datum of finite Cartan type as in the
previous Lemma. Then
\begin{enumerate}
\item $\setst{x_{\beta_1}^{a_1}x_{\beta_2}^{a_2}\ldots x_{\beta_p}^{a_p}}
{a_1, a_2,\ldots , a_p\ge 0}$
forms a basis of $R(\mathcal D)$,
\item $K(\mathcal D)$ is a braided Hopf subalgebra of $R(\mathcal D)$ with
basis
$$\setst{x_{\beta_1}^{Na_1}x_{\beta_2}^{Na_2}\ldots x_{\beta_p}^{Na_p}}{
a_1, a_2, \ldots , a_p\ge 0},$$
\item $[x_{\alpha},x_{\beta}^N]_c=0$, i.e: $x_{\alpha}x_{\beta}^N=
q_{\alpha\beta}^Nx_{\beta}^Nx_{\alpha}$ for all
$\alpha ,\beta \in \Phi^+$.
\end{enumerate}
\end{Theorem}

The vector space $V=V(\mathcal D)$ can also be viewed as a crossed
module in $^{\mathbf Z[I]}_{\mathbf Z[I]}YD$. The Hopf algebra
$\mathcal A(V)$, the quotient Hopf algebra $R(\mathcal D)=\mathcal
A(V)/(ad^{1-a_{ij}}x_i(x_j)|1\le i\ne j\le\theta )$ and its Hopf
subalgebra $K(\mathcal D)$ generated by $\setst{
x_{\alpha}^N}{\alpha\in\Phi^+}$ are all Hopf algebras in
$^{\mathbf Z[I]}_{\mathbf Z[I]}YD$. In particular, their
comultiplications are $\mathbf Z[I]$-graded. By construction, for
$\alpha\in\Phi^+$, the root vector $x_{\alpha} \in R(\mathcal D)$
is $\mathbf Z[I]$-homogeneous of $\mathbf Z[I]$-degree $\alpha$,
so that $x_{\alpha}\in R(\mathcal
D)_{g_{\alpha}}^{\chi_{\alpha}}$. To simplify notation write for
$1\le l\le p$ and for $a =(a_1, a_2,\ldots ,
a_p)\in\mathbf N^p$
$$h_l=g_{\beta_l}^N\ ,\ \eta_l=\chi_{\beta_l}^N\ ,\ z_l=x_{\beta_l}^N$$
and $\underline a=\sum_{i=1}^pa_i\beta_i$
$$h^a =h_1^{a_1}h_2^{a_2}\ldots h_p^{a_p}\in G\ ,\
\eta^a =\eta_1^{a_1}\eta_2^{a_2}\ldots\eta_p^{a_p}\in\tilde G\ , \
z^a =z_1^{a_1}z_2^{a_2}\ldots z_p^{a_p}\in K(\mathcal
D).$$ In particular, for $e_l=(\delta_{kl})_{1\le k\le p}$, where
$\delta_{kl}$ is the Kronecker symbol, $\underline e_l=\beta_l$ and
$z^{e_l}=z_l$ for $1\le l\le p$. The height of $\alpha
=\sum_{i=1}^{\theta}n_i\alpha_i\in\mathbf Z[I]$ is defined to be
the integer $ht(\alpha )=\sum_{i=1}^{\theta}n_i$. Observe that if
$a, b, c\in\mathbf N^p$ and
$\underline a=\underline b+\underline c$ then
$$h^a=h^bh^c\ ,\ \eta^a =
\eta^b\eta^c\ \rm{and}\ ht(\underline
b)<ht(\underline a)\ \rm{if}\ \underline c\ne 0.$$ The
comultiplication on $K(\mathcal D)$ is $\mathbf Z[I]$-graded, so
that
$$\Delta_{K(\mathcal D)}(z^a)=z^a\ot 1 +
1\ot z^a +
\sum_{b, c\neq 0; \underline b+\underline
c=\underline a}t^a_{bc}z^b\ot z^c$$ and hence
$$\Delta_{K(\mathcal D)\#kG}(z^a)=z^a\ot 1 +
h^a\ot z^a +
\sum_{b, c\neq 0; \underline b+\underline
c=\underline a}t^a_{bc}z^bh^c\ot z^c$$ on the
bosonization. The algebra $K(\mathcal D)$ is generated by the
subspace $L(\mathcal D)$ with basis $\set{z_1, z_2,\ldots , z_p}$.
The (left) $kG$-module structure on $\mathcal A(V)$ restricts to
$L(\mathcal D)$, and induces (right) $kG$-actions on $\Alg
(K(\mathcal D),k)$ and on $\Vect (L(\mathcal D),k)$ by the formula
$(fg)(x)=f(gx)$. A linear functional $f\colon L(\mathcal D)\to k$
is called $g$-invariant if $fg=f$ for all $g\in G$. Let
$\Vect_G(L(\mathcal D),k)$ be the subspace of $G$-invariant linear
functionals in $\Vect (L(\mathcal D),k)$.

\begin{Proposition}\label{connected} Let $\Vect_G(L(\mathcal D)$ and
$\Alg_G(K(\mathcal D),k)$ be the space of $G$-invariant linear
functionals and the set of $G$-invariant algebra maps, where
$\mathcal D$ is a connected special datum of finite Cartan type.
Then:
\begin{enumerate}
\item $\Vect_G(L(\mathcal D),k)= \setst{f\in\Vect (L(\mathcal
D),k)}{f(z_l)=0\ \rm{if}\ \eta_l\ne\ep}$.
\item The restriction map $\res\colon \Alg_G(K(\mathcal
D),k)\to\Vect_G(L(\mathcal
D),k)$ is a bijection. The inverse is given by
$\res^{-1}(f)(z^{\underline a})=f(z_1)^{a_1})f(z_2)^{a_2}\ldots
f(z_p)^{a_p}$.
\item $\Alg_G(K(\mathcal D),k)$ is a group under
convolution. \item The restriction map $\res\colon
{_G{\Alg_G}}(K(\mathcal D)\# kG,k)\to \Alg_G(K(\mathcal D),k)$ is
an isomorphism of groups with inverse defined by
$\res^{-1}(f)(x\ot g)=f(x)$, and ${_G\Alg_G}(K(\mathcal D)\#
kG,k)= \setst{\tilde f\in\Alg(K(\mathcal D)\#kG,k)}{\tilde
f_{|kG}=\ep}$.
\item The map $\Theta =\rho\res^{-1}\colon
\Alg_G(K(\mathcal D),k)\to \Aut_{Hopf} (K(\mathcal D)\# kG)^{op}$,
defined by $\Theta (f) =\res^{-1}(f)*1*\res^{-1}(f)s$, is a group
homomorphism whose image is a subgroup in
$$\widetilde{\Aut}_{Hopf}(K(\mathcal D)\# k
G)=\setst{f\in\Aut_{Hopf}(K(\mathcal D)\# k G)}{f_{|kG}=id}.$$
\item For every $f\in\Alg_G(K(\mathcal D),k)$ the automorphism
$\Theta (f)$ of $K(\mathcal D)\#kG$ is determined by
\begin{eqnarray*}
&&\Theta (f)z^a = z^a +f(z^a)(1-h^a)+ \sum_{b, c\ne
0;\underline b+\underline c=\underline a}t^a_{bc}
f(z^b)z^c \\
&& +\sum_{b, c\ne 0;\underline b+\underline
c=\underline a}t^a_{bc}\left[z^b+f(z^b)(1-h^b)
 +\sum_{d, e\ne 0;\underline d+\underline e=\underline b}
t^b_{d,e}f(z^d)z^e\right] h^c fs(z^c).
\end{eqnarray*}
In particular, $\Theta (f)z^a=z^a
+f(z^a)(1-h^a)$ if $ht(\underline a)=1$.
\end{enumerate}
\end{Proposition}

\begin{proof} If $f\in\Vect_G(L(\mathcal D),k)$ then
$f(z_i)=f(gz_i)=\eta_i(g)f(z_i)$ for all
$1\le i\le p$ and for all $g\in G$.
Thus, $f(z_i)=0$ if $\eta_i\ne\ep$.

By Theorem \ref{ASth2.6} it follows that
$$K(\mathcal D)\cong TL(\mathcal D)/(z_iz_j-\eta_j(h_i)z_jz_i|1\le i<j\le
p).$$
If $f\in\Vect_G(L(\mathcal D),k)$ then the induced algebra map
$\tilde f\colon TL(\mathcal D)\to k$ factors uniquely through
$K(\mathcal D)$, since
$$\tilde f(z_iz_j-\eta_j(h_i)z_jz_i)=f(z_i)f(z_j)-
\eta_j(h_i)f(z_j)f(z_i)=f(z_i)(f(z_j)-f(h_iz_j))=0$$
for $1\le i,j\le p$, by the fact that $f$ is $G$-invariant. This
proves the second assertion.

The next three assertions are a special case of \ref{invariant}.

The set of all algebra maps $\Alg (K(\mathcal D),k)$ may not be a
group under convolution, but the subset $\Alg_G(K(\mathcal D),k)$
is. If $f_1$, $f_2$ and $f$ are $G$-invariant then
\begin{eqnarray*}
f_1*f_2(x y) & = & (f_1\ot f_2)(m\ot m)(1\ot c\ot 1)(x_1\ot x_2\ot y_1\ot
y_2) \\
& & =(f_1\ot f_2)(x_1(x_2)_{-1}y_1\ot (x_2)_0y_2 \\
& & =f_1(x_1)\ep ((x_2)_{-1})f_1(y_1)f_2((x_2)_0)f_2(y_2) \\
& & =f_1(x_1)f_1(y_1)f_2(x_2)f_2(y_2) =f_1*f_2(x)f_1*f_2(y)
\end{eqnarray*}
and moreover, $(f_1*f_2)g= f_1g*f_2g =f_1*f_2$, $fsg=fgs=fs$, $\ep
*f =f =f*\ep$ and $f*fs =\ep = fs*f$ so that $\Alg_G(K(\mathcal
D)$ is closed under convolution multiplication and inversion.

The map $\Psi\colon \Alg_G(K(\mathcal D),k)\to \Alg (K(\mathcal
D)\# kG,k)$ given by $\Psi (f)(x\ot g)=f(x)$, is a homomorphism,
since
\begin{eqnarray*}
\Psi (f_1)*\Psi (f_2)(x\ot g)\hskip -2pt &=&
\Psi (f_1)\ot\Psi (f_2)(x_1\ot (x_2)_{-1}g\ot (x_2)_0\ot g) \\
& = & f_1(x_1)\ep ((x_2)_{-1})f_2((x_2)_0) \\
& = & f_1(x_1)f_2(x_2) = f_1*f_2(x) \\
& = & \Psi (f_1*f_2)(x\ot g).
\end{eqnarray*}
The inverse $\Psi^{-1}\colon {_G\Alg_G}(K(\mathcal D)\#
kG,k)\to \Alg_G(K(\mathcal D),k)$, given by $\Psi^{-1}(\tilde
f)(x)=\tilde f(x\ot 1)$, is just the restriction map.

It is convenient to use the notation $\Psi (f)=\tilde f$. Then $\Theta (f)
= \tilde f*1*\tilde fs$ and

\begin{eqnarray*}
\Theta (f_1*f_2) & = & \widetilde{f_1*f_2}*1*\widetilde{f_1*f_2}s \\
& = & (\tilde f_1*\tilde f_2)*1*(\tilde f_2s*\tilde f_1s) \\
& = & \Theta (f_2)\Theta (f_1).
\end{eqnarray*}
In particular,
$\Theta (f)\Theta (fs) =\Theta (fs*f)=
\Theta (\ep )=1=\Theta (f*fs)=\Theta (fs)\Theta (f)$.
Moreover,
\begin{eqnarray*}
\Theta (f) (x y) & = & \tilde f(x_1y_1)x_2y_2\tilde fs(x_3y_3) \\
& = & \tilde f(x_1)\tilde f(y_1)x_2y_2\tilde fs(y_3)\tilde fs(x_3) \\
& = & \tilde f(x_1)x_2\tilde fs(x_3)\tilde f(y_1)y_2\tilde fs(y_3) \\
& =  & \Theta (f)(x)\Theta (f)(y)
\end{eqnarray*}
and
\begin{eqnarray*}
\Delta\Theta (f) & = & \Delta (\tilde f*1*\tilde fs) \\
& = & \Delta (\tilde f\ot 1\ot\tilde fs)\Delta^{(2)} \\
& = & (\tilde f\ot 1\ot 1\ot \tilde fs)\Delta^{(3)} \\
& = & (\tilde f\ot 1\ot\ep\ot 1\ot\tilde fs)\Delta^{(4)} \\
& = & (\tilde f\ot 1\ot\tilde fs\ot\tilde f\ot 1\ot \tilde fs)\Delta^{(5)}
\\
& = & (\tilde f*1*\tilde fs\ot\tilde f*1*\tilde fs)\Delta \\
& = & (\Theta (f)\ot\Theta (f))\Delta ,
\end{eqnarray*}
showing that $\Theta (f)$ is an automorphism of $K(\mathcal D)\#
kG$ with inverse $\Theta (fs)$.

The remaining item now follows from the formula for the
comultiplication
$$\Delta (z^a) = z^a\ot 1 + h^a\ot z^a +
\sum_{b, c\ne 0;\underline b+\underline
c=\underline a}t^a_{bc}z^bh^c\ot z^c$$
of $K(\mathcal D)\# kG$, which implies
\begin{eqnarray*}
&&\hskip -10pt \Delta^{(2)}(z^a) = z^a\ot 1\ot 1 +h^a\ot z^a\ot 1
+ \sum t^a_{bc}z^bh^c\ot z^c\ot 1 + h^a\ot h^a\ot z^a \\
& & + \sum t^a_{bc}\left[z^bh^c\ot h^c\ot
z^c + h^a\ot z^bh^c\ot z^c +\sum t^b_{rl}z^rh^{l+c}\ot z^lh^c\ot
z^c\right].
\end{eqnarray*}
and
$$1*s(z^a)=z^a + h^as(z^a) + \sum t^a_{bc}z^bh^cs(z^c) =\ep (z^a)=0.$$
Applying $\Psi (f)=\tilde f$ to the latter gives
$$f(z^a)+fs(z^a) + \sum t^a_{bc} f(z^b)fs(z^c) = 0,$$
which will be used in the following evaluation. Now compute
\begin{eqnarray*}
\Theta (f)(z^a) & = & f(z^a) + z^a + \sum t^a_{bc}f(z^b)z^c + h^afs(z^a) \\
& & + \sum t^a_{bc}[f(z^b) + z^b + \sum t^b_{rl}f(z^r)z^l] h^cfs(z^c) \\
& = & z^a +f(z^a)(1-h^a) + \sum t^a_{bc}f(z^b) z^c \\
& & + \sum t^a_{bc}[z^b + f(z^b)(1-h^b) +
\sum t^b_{rl}f(z^r)z^l]h^cfs(z^c)
\end{eqnarray*}
to get the required result.
\end{proof}

For any $f\in\Alg_G(K(\mathcal D),k)$ define by induction on
$ht(\underline a)$ the following elements in the augmentation ideal of
$kG$
$$u_a(f)= f(z^a)(1-h^a)
+ \sum_{b, c\ne 0;\underline b+\underline
c=\underline a} t^a_{bc}f(z^b)u_c(f),$$ where $u_a(f)=f(z^a)(1-h^a)$ if
$ht (\underline a)=1$. In particular, for a positive root $\alpha
=\beta_l\in
\Phi^+$ and $x_{\alpha}^N=x_{\beta_l}^N=z^{\underline e_l}$ write
$u_l(f)=u_{\underline e_l}(f)=u_{\alpha}(f)$. We can think of
$f=(f(x_{\alpha}^N)|\alpha\in\Phi^+)$ as root vector parameters in
the sense of \cite{AS}.

\begin{Corollary} Let $\mathcal D$ be a special connected
datum of finite Cartan type. Then
$$u(\mathcal D, f) = R(\mathcal D)\#kG/(x_{\alpha}^N+u_{\alpha}(f))$$
are the liftings of $\mathcal B(V)\#kG =u(\mathcal D, \ep )$.
\end{Corollary}

\begin{proof} The augmentation ideal of $K(\mathcal D)$,
the ideal $I$ of $K(\mathcal D)\#kG$ and the ideal $(I)$ in
$R(\mathcal D)\#kG$ generated by $\setst{x_{\alpha}^N}{\alpha\in
\mathcal X}$ are Hopf ideals. It follows from the inductive
formulas for $\Theta (f)(z^a)$ and $u_a(f)$ above that for every
$f\in\Alg_G(K(\mathcal D),k)$ the
ideals $I_f=\Theta (f)(I)$ in $K(\mathcal D)\#kG$ and $(I_f)$ in
$R(\mathcal D)\#kG$ generated by
$\setst{x_{\alpha}^N+u_{\alpha}(f)}{\alpha\in \Phi^+}$ are Hopf
ideals as well. The Hopf algebras $u(\mathcal D,f)=R(\mathcal
D)\#kG/(\Theta (f)(I))$ are the liftings of $u(\mathcal D,\ep
)=\mathcal B(V)\# kG$ parameterized by
$f=(f(x_{\alpha}^N|\alpha\in\Phi^+)\in\Alg_G(K(\mathcal D),k).$
\end{proof}

In the not necessarily connected case of a special datum of finite
Cartan type the elements $ad^{1-a_{ij}}x_i(x_j)$ are still
primitive in $\mathcal A(V)$ and $R(\mathcal D)=\mathcal
A(V)/(ad^{1-a_{ij}}x_i(x_j)|i\sim j)$ is still a Hopf algebra,
which contains $R(\mathcal D_J)$ for every connected component
$J\in\mathcal X$. The Hopf subalgebra $K(\mathcal D)$ generated by
the subspace with basis  $S=\setst{z_J^{a_J}, z_{ij}}{J\in\mathcal
X, i\not\sim j}$, where $z_{ij}=[x_i,x_j]_c$, contains $K(\mathcal
D_J)$ for every $J\in\mathcal X$. The comultiplication in the
components $K(\mathcal D_J)$ and $K(\mathcal D_J)\#kG$ is of
course given as before in the connected case, while for $i\not\sim j$
$$\Delta (z_{ij})=z_{ij}\ot 1+1\ot z_{ij}$$
in $K(\mathcal D)$ and $R(\mathcal D)$ and
$$\Delta (z_{ij})=z_{ij}\ot 1 + g_ig_j\ot z_{ij}$$
in the bozonizations $K(\mathcal D)\#kG$ and $R(\mathcal D)\#kG$.
The space of $G$-invariant linear functionals $\Vect_G(L(\mathcal
D),k)$ consists of the elements $f\in\Vect (L(\mathcal D),k)$ such that
$$f(z_r)=0\ {\rm if}\ \eta_r\ne\ep\ {\rm for}\ 1\le r\le p\ {\rm and}\
f(z_{ij})=0\ {\rm if}\
\chi_i\chi_j\ne\ep\ {\rm if}\ i\not\sim j.$$
The induced algebra map $\tilde f\colon TL(\mathcal D)\to k$ of
such a linear functional satisfies
\begin{itemize}
\item $\tilde f([z_r,z_s]_c)= f(z_r)(f(z_s)-f(h_rz_s))=0,$
\item $\tilde f([z_{ij},z_r]_c)=f(z_{ij})(f(z_r)-f(g_ig_jz_r)=0,$
\item $\tilde
f([z_{ij},z_{lm}]_c)=f(z_{ij})(f(z_{lm})-f(g_ig_jz_{lm}))=0,$
\end{itemize}
since $f$ is $G$-invariant. It therefore factors through
$K(\mathcal D)$, since
$$TL(\mathcal D)/([z_r,z_s]_c, [z_{ij},z_r]_c, [z_{ij},z_{lm}]_c)
=K(\mathcal D)/([z_r,z_s]_c, [z_{ij},z_r]_c, [z_{ij},z_{lm}]_c).$$
It follows that the restriction maps
$$\res \colon  _G\Alg_G(K(\mathcal D)\#kG,k)\to \Alg_G(K(\mathcal D),k)\to
\Vect_G(L(\mathcal D),k)$$
are bijective, and $f=\setst{
f(z_{ij})}{i\not\sim j}\cup\setst{f(z_r)}{1\le r\le p}$ can be
interpreted as a combination of linking parameters and root vector
parameters in the sense of \cite{AS}. Then map
$$\Theta\colon \Alg_G(K(\mathcal D),k)\to \Aut_{Hopf}(K(\mathcal
D)\#kG)^{op}$$
given by $\Theta (f)=f*1*fs$ is a homomorphism of groups.
Moreover, since $z_{ij}+g_ig_js(z_{ij})=m(1\ot s)\Delta
(z_{ij})=0$ in $K(\mathcal D)\#kG$, it follows that
$$\Theta (f)(z_{ij})=
(f\ot 1\ot fs)\Delta^{(2)}(z_{ij})=z_{ij}+f(z_{ij})(1-g_jg_j)$$
when $i\not\sim j$, while $\Theta f(z_r)$ is given inductively as
in \ref{connected}. In this way one obtains therefore all the
`liftings' of $B(V)\# kG$ for special data of finite Cartan type.

\begin{Theorem}\label{main3} Let $\mathcal D$ be a special datum of
finite Cartan type. Then
$$u(\mathcal D,f)= R(\mathcal
D)\#kG/(x_{\alpha}^{N_{\alpha}}+u_{\alpha}(f),
[x_i,x_j]_c+f(z_{ij})(1-g_ig_j)|\alpha\in\Phi^+ , i\not\sim j)$$
for $f\in\Vect_G(L(\mathcal D),k)$ are the liftings of $\mathcal
B(V)\#kG =u(\mathcal D, \ep )$. Moreover, all these liftings are
monoidally Morita-Takeuchi equivalent.
\end{Theorem}

\begin{proof} Clearly, $u(\mathcal D,f)$ is a lifting of
$\mathcal B(V)\#kG$ for the root vector parameters $\setst{
\mu_{\alpha}=f(x_{\alpha}^{N_{\alpha}})}{\alpha\in\Phi^+}$ and the
linking parameters $\setst{\lambda_{ij}=f([x_i,x_j]_c)}{i\not\sim
j}$. By \cite{AS} all liftings of $\mathcal B(V)\#kG$ are of that
form. To proof the last assertion use \ref{po} with $H=kG$,
$K=K{\mathcal D}$, $f\in\Alg_G(K,k)$. Then the ideal
$I=(x_{\alpha}^{n_\alpha}, [x_i,x_j]_c|\alpha\in\Phi^+, i\not\sim
j)$  and $J=\Theta (f)(I)=(x_{\alpha}^{N_{\alpha}}+u_{\alpha}(f),
[x_i,x_j]_c+f(z_{ij})(1-g_ig_j)|\alpha\in\Phi^+ , i\not\sim j)$ of
$K\#kG$ are conjugate. By \ref{po} the quotient Hopf algebras
$u(\mathcal D,\ep)$ and $u(\mathcal D, f)$ of $R(\mathcal D)\# kG$
are monoidally Morita-Takeuchi equivalent. The additional
condition $(R\# kG)/(\res^{-1}(f)*(I\# kG))\ne 0$ is verified in
(\cite{Ma1}, Appendix).
\end{proof}

\section{Main Results}\label{s3}

In this section we describe liftings of special crossed modules
$V$ over finite abelian groups in terms of cocycle deformations of
$B(V)\#kG$.


A normalized 2-cocycle
$\sigma\colon A\ot A\to k$ on a Hopf algebra $A$ is a convolution
invertible linear map such that
$$(\ep\ot\sigma )*\sigma (1\ot m)=(\sigma\ot\ep )*\sigma (m\ot 1)$$
and $\sigma (\iota\ot 1)=\ep =\sigma (1\ot\iota )$. The deformed
multiplication
$$m_{\sigma}=\sigma *m*\sigma^{-1}\colon A\ot A\to A$$
and antipode
$$s_{\sigma}=\sigma*s*\sigma^{-1}\colon A\to A$$
on $A$, together with the original unit, counit and
comultiplication define a new Hopf algebra structure on H which we
denote by $A_{\sigma}$.

\begin{Theorem} \cite[Corollary 5.9]{Sch} If two Hopf algebras
$A$ and $A'$ are cocycle deformations of each other,
then they are monoidally Morita-Takeuchi equivalent.
The converse is true if $A$ and $A'$ are finite dimensional.
\end{Theorem}

Suppose now that $V$ is a crossed $kG$-module of special finite
Cartan type, $\cal A(V)$ the free braided algebra and $\cal A(V)\#
kG$ its bosonization. If $I$ is the ideal of $\cal A(V)$ generated
by the subset
$$S=\setst{ad^{1-a_{ij}}x_i(x_j)}{i\sim j}\cup
\setst{x_{\alpha}^{N_{\alpha}}}{\alpha\in\Phi^+}\cup\setst{
[x_i,x_j]_c}{i\not\sim j}$$ then $\cal A(V)/I=\cal B(V)$ is the
Nichols algebra. The subalgebra $K$ of $\cal A(V)$ generated by
$S$ is a Hopf subalgebra \cite{AS}, \cite[Proposition 9.2.1]{CP}.
Then $K\# kG$ is the Hopf subalgebra of $\cal A(V)\# kG$ generated
by $S$ and $G$.

\begin{Lemma} The injective group homomorphism
$$\phi \colon \Alg_G(K,k)\to \Alg (K\#kG,k)$$ given by $\phi (f)(x\#
g)=f(x)$ has image $${_G\Alg_G}(K\#kG,k)=\setst{f\in
\Alg(K\#kG,k)}{f_{|kG}=\ep}$$ and
$$\adj\colon \Alg_G(K,k)\to \Aut(K\#kG)^{op}$$
has its image in the subgroup $$\widetilde{Aut(K\#kG)}=\setst{
f\in\Aut (K\#kG)}{f_{|kG}=\ep }.$$ Moreover, if $V$ is of special
finite Cartan type then $f(ad^{1-a_{ij}}x_i(x_j))=0$ for $i\sim j$
and for every $f\in\Alg_G(K,k)$.
\end{Lemma}

\begin{proof} If $f\in \Alg (K\#kG,k)$ then
$f(ad^{1-a_{ij}}x_i(x_j))= f(g\cdot
ad^{1-a_{ij}}x_i(x_j)g^{-1})=\chi_i(g)^{1-a_{ij}}\chi_j(g)$,
$f(x_{\alpha}^N)=f(gx_{\alpha}^{N_{\alpha}}g^{-1})
=\chi_{\alpha}^{N_{\alpha}}(g)f(x_{\alpha}^N)$ and
$f([x_i,x_j]_c)=f(g[x_i,x_j]_cg^{-1})=\chi_i(g)\chi_j(g)f([x_i,x_j]_c)$,
so that $f(g\cdot ad^{1-a_{ij}}x_i(x_j))=0$ if
$\chi_i^{1-a_{ij}}\chi_j\ne 0$, $f(x_{\alpha}^{N_{\alpha}})=0$ if
$\chi_{\alpha}^{N_{\alpha}}\ne\ep $ and $f([x_i,x_j]_c)=0$ if
$\chi_i\chi_j\ne\ep$.
\end{proof}

The theorem above can now be applied to the situation in Section
\ref{s3} to show that all `liftings'  of a crossed $kG$-module of
special finite Cartan type are cocycle deformations of each other.
The special case of quantum linear spaces has been studied by
Masuoka \cite{Ma}, and that of a crossed $kG$-module corresponding
to a finite number of copies of type $A_n$ by Didt \cite{Di}.

\begin{Theorem} \label{main4} Let $G$ be a finite abelian group, $V$ a
crossed
$kG$-module of special finite Cartan type, $\mathcal B(V)$
its Nichols algebra with bosonization $A=\mathcal B(V)\#kG$. Then:
\begin{enumerate}
\item All liftings of $A$ are monoidally Morita-Takeuchi equivalent,
i.e: their comodule categories are monoidally equivalent, or
equivalently,
\item all liftings of $A$ are cocycle deformations of each other.
\end{enumerate}
\end{Theorem}

\begin{proof} Theorem \ref{main3} at the end the last section
says that
$\mathcal B(V)\#kG\cong u(\mathcal D, \ep )\cong R(\mathcal
D)\#kG/(I)$ for a Hopf ideal $I$ in the Hopf subalgebra
$K(\mathcal D)\#kG$ of $R(\mathcal D)\#kG$, that its liftings are
of the form $u(\mathcal D,f)\cong R(\mathcal D)\#kG/(I_f)$ for a
conjugate Hopf ideal $I_f$, where $f\in\Alg_G(K(\mathcal D,k),
k)\cong {_G\Alg_G}(K(\mathcal D)\#kG,k)$, and that they are
all Morita-Takeuchi equivalent. Thus, Schauenburg's result
applies, so that all these liftings are cocycle deformations of
each other.
\end{proof}

\begin{Corollary}\label{main5} Let $H$ be a finite dimensional pointed Hopf algebra
with abelian group of points $G(H)=G$ and assume that the order of $G$ has
no
prime divisors $< 11$. Then:
\begin{itemize}
\item $H$ and $\gr_c(H)$ are Morita Takeuchi equivalent, or equivalently,
\item $H$ is a cocycle deformation of $\gr_c(H)$.
\item There are only finitely many quasi-isomorphism classes of such Hopf
algebras in any given dimension.
\end{itemize}
\end{Corollary}

\begin{proof} Under the present assumptions the Classification Theorem
\cite{AS}
asserts that $\gr_c(H)\cong B(V)\# kG$ for a crossed $kG$-module
$V$ of special finite Cartan type, and hence the Theorem
\ref{main4} applies. Moreover, there are only finitely many isomorphism
classes of Hopf algebras of the form $B(V)\# kG$ in any given dimension.
\end{proof}

\section{About the deforming cocycles}\label{s4}

\subsection{Deforming cocycles and their infinitesimal parts}

Let $\displaystyle{A=\bigoplus_{n=1}^\infty A_n}$ be a $\mathbb{N}$-graded Hopf algebra (here $\mathbb{N}=\{0,1,2\ldots\}$, and $A_n$ denotes the homogeneous component of $A$ of degree $n$).
Let $\sigma\colon A\otimes A\to k$ be a normalized $2$-cocycle as discussed at the beginning of Section \ref{s3}.
\textsl{Throughout the Section we additionally assume that $$\sigma|_{A_0\otimes A_0}=\varepsilon.$$}

We decompose $\sigma=\sum_{n=0}^\infty \sigma_i$ into the (locally finite) sum of homogeneous maps $\sigma_i$ od degree $-i$; more precisely
$$
\sigma_i = A\otimes A \to  \bigoplus_{p+q=i} A_p\otimes A_q = (A\otimes A)_i \stackrel{\sigma|_{(A\otimes A)_i}}{\to} k.
$$
Note that due to our assumption $\sigma|_{A_0\otimes A_0}=\varepsilon$ we have $$\sigma_0=\varepsilon.$$
In this fashion we also decompose $\sigma^{-1}=\sum_{j=0}^\infty \eta_j$, that is $\eta_j$'s are homogeneous
maps uniquely determined by $\displaystyle{\sum_{i+j=\ell} \sigma_i\eta_j = \delta_{\ell, 0}}$ for $\ell=1,2\ldots$
and automatically also satisfy $\displaystyle{\sum_{i+j=\ell} \eta_i\sigma_j = \delta_{\ell, 0}}$.  Note that $\eta_0=\varepsilon$ and that for the least positive integer $s$ for which $\sigma_s\ne 0$ we have $\eta_s=-\sigma_s$.

The cocycle condition
$$(\ep\ot\sigma(t))*\sigma(t)(1\ot m)=(\sigma(t)\ot\ep
)*\sigma(t)(m\ot 1)$$ implies that
$$\sum_{i+j=\ell}(\ep\ot\sigma_i )*\sigma_j(1\ot m)=
\sum_{i+j=\ell}(\sigma_i\ot\ep )*\sigma_j(m\ot 1)$$
for all $\ell\ge 1$. In particular, if $s$ is the least positive
integer for which $\sigma_s\ne 0$, then
$$\ep\ot\sigma_s + \sigma_s (1\ot m)=\sigma_s\ot\ep +\sigma_s(m\ot 1)$$
so that $$\sigma_s\colon A\ot A\to k$$ is a \textsl{Hochschild 2-cocycle}.  We call this Hochschild $2$-cocycle $\sigma_s$ the \textsl{graded infinitesimal part} of $\sigma$.

\subsection{Relationship with formal graded deformations}

The cocycle deformation $A_\sigma$ is a filtered bialgebra with the underlying filtration inherited from the grading on $A$, i.e., the $\ell$-th filtered part is $(A_\sigma)_{(\ell)}=\bigoplus_{i=0}^\ell A_i$.  Note that the associated graded bialgebra $Gr A_\sigma$ can be identified with $A$.  It was observed by Du,Chen, and Ye \cite{DCY} that decomposing the multiplication $m_\sigma=\sigma*m*\sigma^{-1}$ into the sum of homogeneous components $m_i$ of degree $-i$, allows us to identify the filtered $k$-linear structure $A_\sigma$ with a $k[t]$-linear structure $m_\sigma^t \colon A[t]\otimes A[t] \to A[t]$ induced by $m_\sigma^t|_{A\otimes A}=
m+m_1 t +m_2 t^2+\ldots$ (here we assume that the degree of $t$ is 1).  If $\Delta^t\colon A[t]\to A[t]\otimes A[t]$ is the $k[t]$-linear map induced by $\Delta^t|_{A} = \Delta$, then the graded Hopf algebra $A[t]_\sigma= (A[t], m_\sigma^t, \Delta^t)$ is a formal graded deformation of $A$ in the sense of Du, Chen, and Ye \cite{DCY} (see also \cite{GS, MW}).  Note that in case $m$ does not commute with the graded infinitesimal part $\sigma_s$ of $\sigma$, then $(m_s,0)$, where $m_s = \sigma_s*m-m*\sigma_s$, is the infinitesimal part of the formal graded deformation
(and is a $2$-cocycle in the graded version of the truncated Gerstenhaber-Schack bialgebra cohomology).

\subsection{Exponential Map}

It is in general very hard to give explicit examples of multiplicative cocycles. One somewhat
accessible family consists of bicharacters. Below we give another
idea which can sometimes be used.

Note that if $A=\oplus_{n=0}^\infty A_n$ is a graded bialgebra, and
$f\colon A\to k$ is a linear map such that $f|_{A_0}=0$, then
$$
e^f=\sum_{i=0}^\infty \frac{f^{* i}}{i!}\colon A\to k
$$
is a well defined convolution invertible map with convolution
inverse $e^{-f}$.  When $f\colon A\otimes A\to k$ is a Hochschild
cocycle (more precisely, we have $\ep(a)f(b,c)+f(a,bc) = f(a,b)\ep(c)+f(ab,c)$ and
for all $a,b,c\in A$) such that $f|_{A\ot A_0+A_0\ot A}=0$, then `often'
$e^f\colon A\ot A\to k$ will be a multiplicative cocycle. For
instance this happens whenever $f(1\ot m)$ and $f(m\ot 1)$ commute
(with respect to the convolution product) with $\ep\ot f$ and
$f\ot \ep$, respectively. Also note that if $f*f=0$, then
$e^f=\ep+f$.

\begin{Lemma}
If $f\colon A\otimes A\to k$ is a Hochschild $2$-cocycle such that $f(1\ot m)$ commutes with $\ep\ot f$ and
$f(m\ot 1)$ commutes with $f\ot\ep$ in the convolution algebra $\Hom_k(A\otimes A\otimes A, k)$, then
$e^f$ is a (multiplicative) $2$-cocycle with graded infinitesimal part equal to $f$.
\end{Lemma}
\begin{proof}
Since $f$ is a Hochschild $2$-cocycle we have $\ep\ot f + f(1\ot m) = f\ot \ep + f(m\ot 1)$.  Note that
$e^{\ep\ot f} = \ep \ot e^f$ and $e^{f\ot\ep} = e^f\ot \ep$. Since $m$ is a coalgebra map we also have $e^{f(m\ot 1 )}
= e^f(m\ot 1)$ and $e^{f(1\ot m )} = e^f(1\ot m)$.  Also recall that for commuting $g and h$ we have $e^{g+h}=e^g* e^h$.
Hence we have
\begin{eqnarray*}
(\ep\ot e^f)*\left(e^f(1\ot m)\right)&=& e^{\ep\ot f}* e^{f(1\ot m)} \\
&=& e^{\ep\ot f+ f(1\ot m)} \\
&=& e^{f\ot \ep+ f(m\ot 1)} \\
&=& e^{f\ot \ep}* e^{f(m\ot 1)} \\
&=& (e^f\ot \ep)*\left(e^f(m\ot 1)\right).
\end{eqnarray*}
\end{proof}

\subsection{Exponential map in a quantum linear space}
In the rest we will study $u(\D,f)$ where the Dynkin diagram
associated to $\D$ is of type $A_1\times\ldots \times A_1$.  We study the linking and root-vector parts of $f$ separately
and use notation $u(\D,f)=u(\D,\lambda,\mu)$, where $\mu=(\mu_\alpha)_{\alpha\in\Phi^+} = (f(x_\alpha^{N_\alpha}))_{\alpha\in\Phi^+}$ (root-vector parameters), and $\lambda=(\lambda_{i,j})_{1\le i<j<\theta} = (f([x_i,x_j]_c))_{1\le i<j\le\theta}$ (linking parameters).

From now on assume that $A=u(\D,0,0)$ is obtained as the bosonization of a quantum
linear space. More precisely,
$$A= B(V)\# kG=\genst{G, x_1,\ldots , x_\theta}{gx_i =\chi_i(g) x_i g,
x_ix_j = \chi_j(g_i) x_j x_i, x_i^{N_i}=0}.$$
Here $V=\oplus_{i=1}^{\theta}kx_i$, and $\chi_1,\ldots ,
\chi_\theta\in\widehat G$, $g_1,\ldots , g_\theta\in G$
are such that $\chi_i(g_j)\chi_j(g_i)=1$ for $i\not=j$, and the number
$N_i$ is the order of $\chi_i(g_i)$. Recall that
$q_{ij}=\chi_i(g_j)$, and let $e_i=(\delta_{ij})_{j=1}^{\theta}$ for $1\le
i\le\theta$, where $\delta_{ij}$ is the Kronecker symbol. Then
$\zeta_i\colon A\ot A\to k$, given by
$$\zeta_i(x^ag, x^bh)=\begin{cases}\chi_i^{b_i}(g) ,& \mbox{ if }
a=a_ie_i, b=b_ie_i, a_i+b_i=N_i \\
0 ,& \mbox{ otherwise }\end{cases}$$
for $x^a=x_1^{a_1}\ldots
x_\theta^{a_\theta}$ and $x^b=x_1^{b_1}\ldots x_\theta^{b_\theta}$
(see for example \cite{MW}) is a Hochschild cocycle.  Moreover, each
of the sets $$S_l=\setst{(\ep\ot\zeta_i), \zeta_i(1\ot m)}{1\le
i\le\theta}$$ and
$$S_r=\setst{(\zeta_i\ot \ep), \zeta_i(m\ot 1)}{1\le
i\le\theta}$$ is a commutative set (for the convolution product). We
sketch the proof for $S_l$ (the proof for $S_r$ is symmetric).
The maps $\zeta_i(1\ot m)$ and $\zeta_j(1\ot m)$ commute since
$\zeta_i$ and $\zeta_j$ do. The same goes for $\ep\ot\zeta_i$ and
$\ep\ot\zeta_j$. Hence it is sufficient to prove that for all
$i,j$ we have
$$
(\ep\ot\zeta_i)*(\zeta_j(1\ot m))=  (\zeta_j(1\ot
m))*(\ep\ot\zeta_i).
$$
If $i\not=j$, this is immediate. For $i=j$ note that both left and
right hand side can be nonzero only at PBW elements of the form
$x_i^r f\ot x_i^s g\ot x_i^p h \in A\ot A\ot A$, with
$r+s+p=2N_i$. Without loss of generality assume that $f=g=h=1$. In
this case the left hand side evaluates to
\begin{eqnarray*}
\sum_{u+v=N_i} {s\choose u}_{q_{ii}} {p\choose v}_{q_{ii}}
q_{ii}^{u(p-v)} = 1
\end{eqnarray*}
and the right hand side is
\begin{eqnarray*}
\sum_{u+v=N_i} {r\choose u}_{q_{ii}} {s\choose v}_{q_{ii}}
q_{ii}^{u(s-v)} = 1.
\end{eqnarray*}
Thus we have the following.

\begin{Proposition}\label{P-root}
If $\zeta = \sum_{i} \mu_i \zeta_i$, (where $\zeta_i$'s are as above),
then $\zeta\colon A\ot A\to k$ is a Hochschild $2$-cocycle, and $\sigma
=e^{\zeta}$ is a multiplicative cocycle. The associated cocycle
deformation is a lifting of $A$ with presentation
\begin{eqnarray*}
\displaystyle
A_{\sigma} &=& u(\D,0,\mu) \\
&=&\genst{G, x_1,\ldots ,
x_\theta}{gx_i =\chi_i(g) x_i g, x_ix_j = \chi_j(g_i) x_j x_i,
x_i^{N_i}=\mu_i(1-g_i^{N_i})}.
\end{eqnarray*}
\end{Proposition}
\begin{proof}
Note that for $1\le i<j\le \theta$ we have $m_{A_\sigma}(x_i\ot x_j - \chi_{j}(g_i) x_j\ot x_i)=0$ and $m_{A_\sigma}(x_i^{\ot N_i})=\mu_i (1-g_i^{N_i})$.  Hence $A$ is a quotient of $u(\D,0,\mu)$.  Since the dimensions of $A_\sigma$ and $u(\D,0,\mu)$ are the same we conclude that
$A_\sigma = u(\D,0,\mu)$.
\end{proof}

\subsection{q-exponential map}

Here we use the ideas and results on q-exponential maps of Sarah Witherspoon \cite{W} in order to address the linking cocycles.
The second author thought of this idea in a conversation with Margaret Beattie about
the papers \cite{GM1} and \cite{ABM}.   The first author obtained a similar idea independently and earlier through cohomological considerations \cite{Gr2}.

If $q$ is a root of $1$ and $x$ is an element of any algebra $B$, then we can define a
$q$-exponential function as follows.
\begin{Definition}
If $\ell\ge 2$ and $q$ is a primitive $\ell$-th root of unity then $\exp_q(x)=\sum_{m=0}^{\ell-1}
\frac{1}{(m)!_q} x^m$.
\end{Definition}

We will need the following nice addition formula proved in \cite{W}.
\begin{Lemma}[\cite{W}]
If  $\ell\ge 2$, $q$ is a primitive $\ell$-th root of unity, $yx=qxy$, and $x^iy^{\ell-i}=0$ for
$i=0,\ldots, \ell$, then $\exp_q(x+y)=\exp_q(x)\exp_q(y)$.
\end{Lemma}

From now on assume, as in the subsection above, that $A$ is a bosonization of a quantum linear space.  Extend the characters $\chi_i$ to all of $A$ by letting
$$\chi_i(xg)=\ep(x)\chi_i(g)$$ for any PBW-element $x\in A$. Define skew-derivations $d_i\colon A\to k$ by setting
$$d_i((x_jg)=\delta_{ij},$$
so that $\chi_im =\chi_i\ot \chi_i$ and $d_im=d_i\ot\ep + \chi_i\ot d_i$. Consider
the linear maps
$$\eta_{j,i}=d_j*\chi_i\ot d_i\colon A\ot A\to k.$$
for which we have the following result.

\begin{Remark}If $i, j$ are linkable roots (i.e., $\chi_i\chi_j=\varepsilon$, $g_ig_j\not=1$, and $\chi_i(g_j)\chi_j(g_i)=1$), then $\eta_{j,i}=d_j*\chi_i\otimes d_i$ is a Hochschild cocycle. More precisely, it is the $2$-cocycle on $A$ associated to the cup products of the Hochschild  $1$-cocycles $d_j$, $d_i$ on $A$; for more details see \cite{GM1, MW, MPSW}.
\end{Remark}

We will now show that $\exp_q(\eta_{j,i})$ is a multiplicative cocycle (the q-exponential is defined with respect to the convolution product in $(A\otimes A)^*=A^*\otimes A^*$).
\begin{Theorem}
If $i,j$ are linkable roots, $q=\chi_i(g_j)$, and $\eta_{j,i}\colon A\otimes A\to k$ is as above,
then $\sigma:=exp_q(\eta_{j,i})\colon A\otimes A\to k$ is a multiplicative cocycle.
\end{Theorem}
\begin{proof} We need to prove that
$$
(\varepsilon\otimes \sigma)*\sigma(\id\otimes\m) = (\sigma\otimes\varepsilon)*\sigma(\m\otimes\id).
$$
First note that $\varepsilon\otimes\sigma =\exp_q(\varepsilon\otimes\eta_{j,i})$, $\sigma\otimes\varepsilon =\exp_q(\eta_{j,i}\otimes\varepsilon)$, $\sigma(\id\ot\m) = \exp_q(\eta_{j,i}(\id\ot\m))$, and
$\sigma(\m\ot\id) = \exp_q(\eta_{j,i}(\m\ot\id))$. The $q$-exponentials here are defined in the convolution algebra $(A\ot A)^*$.
Observe that  for $r\in\{i,j\}$ we have $d_r\m = \Delta_{A^*} d_r = d_r\ot\ep + \chi_r\ot d_r$ and that $\chi_r\m =
\Delta_{A^*} \chi_r = \chi_r\ot \chi_r$.  Hence we have $\varepsilon\otimes\eta_{j,i}=\ep\ot d_j*\chi_i\ot d_i$,
$\eta_{j,i}(\id\ot\m) = d_j*\chi_i\ot d_i\ot\ep + d_j*\chi_i\ot\chi_i\ot d_i$, $\eta_{j,i}\otimes\varepsilon =
d_j*\chi_i\ot d_i\ot \ep$, and $\eta_{j,i}(\m\ot\id) = d_j*\chi_i\ot\chi_i\ot d_i+ \ep\ot d_j*\chi_i\ot d_i$.  If
$a=\ep\ot d_j*\chi_i\ot d_i$, $b=d_j*\chi_i\ot d_i\ot\ep$, $c=d_j*\chi_i\ot\chi_i\ot d_i$ then an easy calculation shows that
$a*b=b*a$, $c*a=qa*c$, $c*b=qb*c$, and $a^{i}*c^{\ell-i}=0=b^i*c^{\ell-i}$ for $i=0,\ldots,\ell$.  By Lemma 4.4 it now follows that
\begin{eqnarray*}
\lefteqn{(\varepsilon\otimes \sigma)*\sigma(\id\otimes\m) = \exp_q(a)*\exp_q(b+c)}\\
&=& \exp_q(a)*\exp_q(b)*\exp_q(c) \\
&=& \exp_q(b)*\exp_q(a)*\exp_q(c) \\
&=& \exp_q(b)*\exp_q(a+c) \\
&=& (\sigma\otimes\varepsilon)*\sigma(\m\otimes\id)
\end{eqnarray*}
as required.
\end{proof}

\begin{Lemma}\label{lxx}
If $(i,j)$ and $(s,t)$ are disjoint pairs of linkable roots, then $\eta_{j,i}$ and $\eta_{s,t}$ commute, and thus also $\exp_{\chi_i(g_j)}(\eta_{j,i})$ and $\exp_{\chi_s(g_t)}(\eta_{t,s})$ commute.
\end{Lemma}
\begin{proof}
It follows directly from the fact that
$$
\chi_i(g_t)\chi_t(g_i)\chi_j(g_s)\chi_s(g_j)=1.
$$
\end{proof}

Let $G_0=\setst{ g_\beta^{N_\beta}}{\beta\in\Phi^+, \chi_\beta^{N_\beta}=\ep, g_\beta^{N_\beta}\not=1}$.  Let $\D'$ be the datum obtained from $\mathcal{D}$ by replacing $G$ by $G'=G/G_0$. It is well know \cite{AS} that elements of $G_0$ are central in $u(\D, \lambda, \mu)$ for all
$\lambda$ and $\mu$.  Abbreviate $I_\mu$ the ideal in $u(\D,0,\mu)$ generated by $(kG_0)^+$ and note that $u(\D,0,\mu)/I_\mu = u(\D',0,0)$. We will use $\pi_\mu$ to denote the canonical projection $\pi_\mu\colon u(\D,0,\mu)\to u(\D,0,\mu)/I_\mu = u(\D',0,0)$.

\begin{Proposition}\label{P-linking} Let $\lambda$ be a linking datum for $\D$ (and hence also for $\D'$).
Let $\sigma_\lambda\colon u(\D',0,0)\otimes u(\D',0,0)\to k$ be the multiplicative cocycle from Lemma \ref{lxx}. Then $\sigma_{\lambda,\mu}:=\sigma_\lambda\circ(\pi_\mu\otimes \pi_\mu)\colon u(\D,0,\mu)\otimes u(\D,0,\mu)\to k$ is a multiplicative cocycle and
$u(\D,0,\mu)_{\sigma_{\lambda,\mu}} = u(\D,\lambda,\mu)$.
\end{Proposition}
\begin{proof}  Proof is almost identical to the proof of Proposition \ref{P-root}.
Let $B=u(\D,0,\mu)_{\sigma_{\lambda,\mu}}$.  Note that for $1\le i<j\le \theta$ we have $m_B(x_i\ot x_j - \chi_{j}(g_i) x_j\ot x_i)=\lambda_{i,j}(1-g_ig_j)$ and $m_B(x_i^{\ot N_i})=\mu_i (1-g_i^{N_i})$.  Hence $B=u(\D,0,\mu)_{\sigma_{\lambda ,\mu}}$ is a quotient of $A=u(\D,\lambda,\mu)$.  Since the dimensions of $A$ and $B$ are the same we conclude that $u(\D,0,\mu)_{\sigma_{\lambda ,\mu }}=u(\D,\lambda,\mu)$.
\end{proof}

We can now explicitly describe all deforming cocycles for the quantum linear space.
\begin{Theorem}\label{thm-explicit}
Let $\sigma:=\sigma_{\lambda,\mu}*e^{\sum_{\beta} \mu_i\zeta_i}$.  Then $\sigma$ is a multiplicative cocycle and $u(\D,0,0)_\sigma =
u(\D,\lambda,\mu)$.
\end{Theorem}
\begin{proof}
Follows from Propositions \ref{P-linking} and \ref{P-root}.
\end{proof}

\subsection{Remarks regarding the general case}\label{general}
The q-exponential map idea described above can be adjusted \cite{GM2} to the general case (liftings of bosonizations of Nichols algebras of finite Cartan type) and one can thus obtain the explicit description of the cocycles deforming $u(\mathcal{D},0,\mu)$ to $u(\mathcal{D},\lambda,\mu)$.  In other words Proposition \ref{P-linking} extends to liftings of bosonizations of all Nichols algebras of finite Cartan type.  The procedure is almost the same as the procedure for a quantum linear space.  The only subtle point is to a carefully order the linkable roots in such a way that all $q$-exponentials involved commute (replacing $\exp_q(d_i\chi_j\ot d_j)$ by $\exp_q(d_j\chi_i\ot  d_i)$ when necessary) .  The same idea also works to explicitly describe the cocycles that deform the `big' quantum groups, that is deforming $U(\mathcal{D},0)$ to $U(\mathcal{D},\lambda)$.

Recently we have also been able to adapt our exponential map ideas to explicitly describe cocycles deforming $u(\mathcal{D},0,0)$ to $u(\mathcal{D},0,\mu)$ for a much larger class of Andruskiewitsch-Schneider Hopf algebras.  We are preparing this work as a separate publication \cite{GM2}.  There we also explain how all of the root vector cocycles can be obtained from Singer cocycles.  We briefly describe this below:
Let $G_0$ be the subgroup of $G$ generated by those $g_i^{N_i}$ for which $\mu_i\not=0$ and let $A=kG_0$.  Note that
$A$ is a central Hopf subalgebra of $B=u(\mathcal{D},0,\mu)$ and that the `quotient' is $C=u(\mathcal{D}',0,0)$.  We thus get
a Singer extension (also called cleft extension) \cite{Si}  $A\stackrel{i}{\to} B\stackrel{\pi}{\to} C$.  Let $a\colon C\otimes C\to A$ be the algebra part of the associated Singer cocycle (it can be computed explicitly; see below).  Then for any character $\phi\colon G_0\to k$ we get a multiplicative cocycle $\sigma'=\phi\circ a\colon C\otimes C\to k$.  This cocycle naturally lifts from $C=u(\mathcal{D}',0,0)$ to $u(\mathcal{D},0,0)$ and the cocycle twist with the above multiplicative cocycle yields $u(\mathcal{D},0,\mu')$, where the datum $\mu'$ is the appropriately re-scaled datum $\mu$ (depending on the choice of the character $\phi$).

How to get $a$ explicitly \cite{An, By}:  Let $\xi\colon B\to A$ denote the section of $i$  that sends PBW-elements from $G_0$ to themselves and is $\varepsilon$ on all other PBW-elements.  As $\xi$ is a unit and counit preserving $A$-module map we can define a $C$-comodule map $\chi\colon C\to B$ by $\chi\pi = (i\circ \xi^{-1})*id$.  Map $a\colon C\otimes C\to A$ is then given by $a(x,y)=\chi(x_1)\chi(y_1)\chi^{-1}(x_2 y_2)$ (the range falls into the $C$-coinvariant part of $B$, which is $i(A)$).

Combining the linking cocycles obtained by the $q$-exponential formula with the root-vector cocycles obtained from Singer cocycles then describes explicitly all deforming cocycles for all of the Andruskiewitsch-Schneider Hopf algebras.

\section{Examples}\label{s5}

The following examples of explicit description of
deforming cocycles of liftings of various quantum linear spaces illustrate the ideas and discussion above.  The first
(and simplest) example is worked out in a lot of detail, implicitly repeating many of the more general constructions above.  This was originally done by the authors in order to better understand the lifting process and gain insight into the structure of liftings.  We hope that the reader will similarly benefit from the experience.
\begin{Example}\label{example1}

{\rm Let $G=\langle g\rangle$ be a cyclic group of order $Np$ and
let $V=kx$ be a 1-dimensional crossed $G$-module with action and
coaction given by $gx=qx$ for a primitive $N$-th root of unity $q$
and $\delta (x)=g\ot x$. The braiding $c\colon V\ot V\to V\ot V$
is then determined by $c(x\ot x)=qx\ot x$. The braided Hopf
algebra $\mathcal A(V)$ is the polynomial algebra $k[x]$ with
comultiplication $\Delta (x^i)=\sum_{r+s=i}{i\choose r}_q x^r\ot
x^s$ in which $x^N$ is primitive. The braided Hopf algebra
$\mathcal C(V)=k\langle x\rangle$ is the divided power Hopf
algebra with basis $\setst{x_i}{i\ge 0}$, comultiplication $\Delta
(x_i)=\sum_{r+s=i}x_r\ot x_s$ and multiplication
$x_ix_j={{i+j}\choose i}_qx_{i+j}$. The quantum symmetrizer
$\mathcal S\colon \mathcal A(V)\to\mathcal C(V)$ is given by
$\mathcal S (x^i)=\mathcal S(x)^i=i_q!x_i$. The Nichols algebra of
$V$ is $B(V)=\mathcal A(V)/(x^N)\cong\im\mathcal S$ and the Hopf
algebra
\begin{eqnarray*}
A &=& B(V)\# k G \\
&=&\left\langle x, g|x^N=0, g^{Np}=1, gx=qxg, \Delta (x)=x\ot
1+g\ot x, \Delta (g)=g\ot g\right\rangle
\end{eqnarray*}
is coradically graded.

The linear functional $\zeta \colon A\ot A\to k$ of degree $-N$,
defined by
$$\zeta (x^ig^j\ot x^kg^l)=aq^{jk}\delta^{i+k}_N,$$
is a Hochschild cocycle with $\zeta^2=0$, and satisfying
$$\zeta (m\ot 1)*(\zeta\ot\ep )=\zeta (1\ot m)*(\ep\ot\zeta ).$$
It follows that
$$\sigma =e^{\zeta}=\ep\ot\ep +\zeta\colon A\ot A\to k$$
is a convolution invertible multiplicative cocycle. In
terms of the dual basis of $A^*$ it can be expressed as
$$\sigma =\ep\ot\ep +\zeta =\ep\ot\ep +
\sum_{r+s=N}^{0<r,s<N}a_{rs}\xi^r\theta^{s}\ot\xi^s ,$$
where $a_{rs}={1\over{r_q!s_q!}}$. The corresponding cocycle
deformation of the multiplication of $A$ is
$$m_{\sigma}=(\sigma\ot m\ot\sigma^{-1})\Delta^{(2)}_{A\ot A}=
m+(\zeta\ot m-m\ot\zeta )\Delta_{A\ot A},$$
since $(\zeta\ot m\ot\zeta )\delta_{A\ot A}=0$ (it is of degree $-2N$).
Using
$$\Delta_{A\ot A}(x^ig^j\ot x^kg^l)=
\sum_{r+s=i}^{u+v=k}{i\choose r}_q{k\choose u}_qx^rg^{s+j} \ot
x^ug^{v+l}\ot x^sg^j\ot x^vg^l$$ and invoking the identity
\cite{Ka}
$$\sum_{s+v=\beta}{i\choose s}_q{k\choose v}_qq^{s(k-v)}=
{{i+k}\choose\beta }_q={{N\alpha +\beta }\choose\beta }_q=1$$ when
$i+k=N\alpha +\beta $, the following explicit formula for $m_{\sigma}$
can be deduced:
\begin{eqnarray*} &&\hskip -10pt
m_{\sigma}(x^ig^j\ot
x^kg^l)=q^{jk}x^{i+k}g^{j+l}\\
&&+\sum_{r+s=i}^{u+v=k}a{i\choose r}_q{k\choose
u}_qq^{(s+j)u+jv} (\delta^{r+u}_{n\alpha}x^{s+v} -
x^{r+u}\delta^{s+v}_{n\alpha}g^{s+v})g^{j+l}\\
&&=q^{jk}\left[x^{i+k}+ax^{\beta} (\sum_{s+v=\beta}{i\choose
s}_q{k\choose v}_qq^{s(k-v)}\right.\\
&&\hskip 10pt - \left.\sum_{r+u=\beta}{i\choose
r}_q{k\choose u}_qq^{(i-r)u}g^{N\alpha}\right]g^{j+l}\\
&&=q^{jk}\left[x^{i+k}+ax^{\beta}(1-g^{N\alpha})\right]g^{j+l},
\end{eqnarray*}
where $i+k=N\alpha +\beta$ with $\alpha =0, 1$.

The linear dual $V^*=k\xi$, $\xi (x)=1$, is a crossed module over
the character group $\widehat{G}=\gen{\theta}$, $\theta
(g)=\alpha$ a primitive $Np$-th root of unity and $\alpha^p=q$,
with action $\delta^*(\theta\ot\xi )=\theta\xi =\alpha\xi$ and
coaction $\mu^*(\xi )=\phi\ot\xi$, where $\phi =\theta^p$. The
graded braided Hopf algebra $\mathcal A(V^*)\cong k[\xi
]\cong\mathcal C(V)^*$ is the graded polynomial algebra with
comultiplication $\Delta (\xi^i )=\sum_{r+s=i}{i\choose r}_q
\xi^r\ot \xi^s$ so that $\xi^N$ is primitive. The cofree graded
braided Hopf algebra $\mathcal C(V^*)=k\gen{\xi} \cong\mathcal
A(V)^*$ is the divided power Hopf algebra with basis
$\setst{\xi_i}{i\ge 0}$, comultiplication $\Delta
(\xi_i)=\sum_{r+s=i}\xi_r\ot\xi_s$ and multiplication
$\xi_i\xi_j={{i+j}\choose i}_q\xi_{i+j}$. The quantum symmetrizer
$\mathcal S\colon \mathcal A(V^*)\to \mathcal C(V^*)$ is given by
$\mathcal S(\xi^i)=i_q!\xi_i$. The Nichols algebra of $V^*$ is
$B(V^*)=\mathcal A(V^*)/(\xi^N)\cong\im\mathcal S$ and the Hopf
algebra
\begin{eqnarray*} A^*&=&B(V^*)\# k\widehat{G} \\
&=& \genst{\xi, \theta}{\xi^N=0, \theta^{Np}=\ep , \theta \xi
=\alpha \xi\theta , \Delta (\xi )=\xi\ot\ep +\phi\ot\xi , \Delta
(\theta )=\theta\ot\theta } \end{eqnarray*}
is radically graded.

The invertible element $\sigma^*\colon k\to A^*\ot A^*$ with
$\sigma^*(1)=\sigma =\ep\ot\ep
+\sum_{r+s=n}a_{rs}\xi^r\phi^s\ot\xi^s =\ep\ot\ep +\zeta$, with
$a_{rs}={1\over{r_q!s_q!}}$, is the cocycle above represented in
terms of the basis of $A^*$. Observe that $\zeta$ is of degree $n$
and $\zeta^2=0$. The resulting cocycle deformation of the
comultiplication
$$\Delta_{\sigma}=m^{(2)}_{A\ot A}\left(\sigma\ot\Delta\ot\sigma^{-1}\right)=
\Delta +m_{A\ot A}(\zeta\ot\Delta -\Delta\ot\zeta
)=\delta_0+\delta_n,$$ where $m^{(2)}_{H\ot
H}(\zeta\ot\Delta\ot\zeta )=0$ is used, is compatible with the
original multiplication. Since $\Delta (\theta )\zeta
=\alpha^N\zeta\Delta (\theta )$, it follows that
$$\Delta_{\sigma}(\theta^i)=\theta^i\ot\theta^i +
(1-\alpha^{Ni})\zeta (\theta^i\ot\theta^i).$$
Using the identity $a_{u-1,s}q^s+a_{u,s-1}=a_{u,s}$ one finds that
$\zeta\Delta (\xi )=\Delta(\xi )\zeta$, so that
$\Delta_{\sigma}(\xi )=\Delta (\xi )$ and
$$\Delta_{\sigma }(\xi^i\theta^j)=\Delta
(\xi^i)\Delta_{\sigma}(\theta^j).$$
The deformed Hopf algebra has the presentation
$$A^{\sigma}=\genst{\xi , \theta}{\Delta_{\sigma}(\xi )=\Delta (\xi ),
\Delta_{\sigma}(\theta )= \theta\ot\theta +(1-\alpha^N)\zeta
(\theta\ot\theta )}$$ with the original multiplication and radical
filtration, so that $\gr_rA^{\sigma }=A^*$. }
\end{Example}

\begin{Example}[cf. \cite{ABM}]\label{E2}{\rm

This is a $2$-variable analogue of the example above. Let
$G=\gen{g}$ be the cyclic group of order $mn$ and $q$ a
primitive $m$-th root of unity. Consider the $2$-dimensional
crossed $G$ module $V=kx_1\oplus kx_2$ with action
$gx_i=q^{(-1)^{i-1}}x_i$ and coaction $\delta (x_i)=g\ot x_i$. The
braiding map $c\colon V\ot V\to V\ot V$ is then $c(x_i\ot
x_j)=q^{(-1)^{j-1}}x_j\ot x_i$.

The bosonization of the Nichols algebra is then
$$A=B(V)\# kG=\genst{g, x_1, x_2}{g^{mn}=1, x_1^m=0, x_2^m=0,
x_2x_1=q x_1x_2}$$
with comultiplication $\Delta (g)=g\ot g$ and $\Delta (x_i)=x_i\ot
1+g\ot x_i$. Note that it is coradically graded.

The linear functionals $\zeta_i\colon A\ot A\to k$ of degree $-m$,
given by $$\zeta_i(x_1^{j_1} x_2^{j_2} g^r\ot
x_1^{j_1'}x_2^{j_2'}g^{r'})= \delta_{j_i+j'_i,m}\delta_{j_{i+1}+
j'_{i+1},0}q^{(-1)^irj'_i}$$ are commuting Hochschild cocycles and the
cocycle deformation associated to the multiplicative cocycle
$\sigma_{\mu_1,\mu_2}=e^{\mu_1\zeta_1+\mu_2\zeta_2} = \ep +\mu_1\zeta_1 + \mu_2\zeta_2 +\mu_1 \mu_2\zeta_1\zeta_2$ has presentation
$$A_{\sigma_{\mu_1,\mu_2}}=\genst{g, x_1, x_2}{g^{mn}=1, x_1^m=\mu_1(1-g^m), x_2^m=
\mu_2(1-g^m), x_2x_1=q x_1 x_2}.$$
Moreover, if $d_i\colon A\to k$ are skew-derivations given by
$$d_1(x_1^{a_1}x_2^{a_2} g) = \begin{cases} 1; a_1 = 1, a_2 = 0 \\ 0 ;\mbox{ otherwise}\end{cases}  , \  d_2(x_1^{a_1}x_2^{a_2} g) = \begin{cases} 1; a_1 = 0, a_2 = 1 \\ 0 ;\mbox{ otherwise}\end{cases}$$
then $\sigma=\sigma_{\lambda; \mu_1,\mu_2} = \exp_q(\lambda d_2\chi\ot d_1)*\sigma_{\mu_1, \mu_2}$ is a multiplicative cocycle
and $A_\sigma$ has presentation
$$\genst{g, x_1, x_2}{g^{mn}=1, x_1^m=\mu_1(1-g^m), x_2^m=
\mu_2(1-g^m), x_2x_1-qx_1x_2 = \lambda(1-g^2)}.$$
}
\end{Example}

\end{document}